\documentclass{article}

\usepackage{amsmath,yhmath,amsfonts,amssymb,amsthm,bbding}
\usepackage{mathrsfs}

\newcommand{\dist}{{\operatorname{dist}}}

\newcommand{\sn}{{\operatorname{sn}}}

\newcommand{\ct}{{\operatorname{ct}}}

\newcommand{\counte}{theorem}
\newtheorem{theorem}{\bf Theorem}[section]

\newtheorem{defn}[\counte]{\bf Definition}

\newtheorem{prop}[\counte]{\bf Proposition}
\newtheorem{lemma}[\counte]{\bf Lemma}

\newtheorem{remark}[\counte]{\bf Remark}

\newtheorem{prop*}{\bf Step} %???????

\allowdisplaybreaks

\numberwithin{equation}{section}

\renewcommand{\thefootnote}{\fnsymbol{footnote}}

\begin{document}

\renewcommand{\thefootnote}{\arabic{footnote}}

\centerline{\bf\Large New definitions of Alexandrov space and applications \footnote{Supported by NSFC 11971057 and BNSF Z190003. \hfill{$\,$}}}

\vskip5mm

\centerline{Shengqi Hu, Xiaole Su,  Yusheng Wang\footnote{The
corresponding author (E-mail: wyusheng@bnu.edu.cn). \hfill{$\,$}}}

\vskip6mm

\noindent{\bf Abstract.} In this paper we show that, in the definition of Alexandrov spaces with lower or upper curvature bound,
the original conditions can be replaced with much weaker ones. For the purpose, we introduce `imaginary' comparison angles (and `imaginary' angles),
and the right or left bounded second derivative in the support sense. As applications, we provide new proofs for the Doubling Theorem,
and the Globalization Theorem for complete or geodesic Alexandrov spaces with lower curvature bound.
\vskip1mm


\vskip1mm


\vskip6mm

\setcounter{section}{-1}

%%%%%%%%%%%%%%%%%%%%%%%%%%%%%%%%%%%%% section 0 Introduction%%%%%%%%%%%%%%%%%%%%%%%%%%%%%%%%%%%%%%%%%

\section{Introduction}

Alexandrov spaces with lower curvature bound have been developed systematically by  Burago-Gromov-Perel$'$man
in \cite{BGP}. Afterwards, Alexandrov geometry has been studied intensively by many geometrists.
With respect to $\mathbb S^2_k$, the complete and simply connected $2$-dimensional space form with constant curvature $k$,
we can define an intrinsic metric space to be an Alexandrov space with curvature $\geqslant k$.
There are several equivalent conditions in \cite{BGP} to define such a kind of spaces.
We first recall the one for general intrinsic metric spaces.

\begin{defn}{\rm
	An intrinsic metric space $X$ is called an {\it Alexandrov space with curvature} $\geqslant k$, if around any $x\in X$ there is a neighborhood $U_x$ such that the following condition is satisfied:
	\vskip1mm
	
	\noindent{\rm (1)}\ \ For any distinct points $a,b,c,d\in U_x$, it holds that
	\begin{align}
		\widetilde{\angle}_kbac+\widetilde{\angle}_kbad+\widetilde{\angle}_kcad\leqslant 2\pi,\label{0-D}
	\end{align}
where $\widetilde\angle_kbac$ denotes the angle at $\tilde a$ of a triangle $\triangle \tilde a\tilde b\tilde c\subset\mathbb S^2_k$
with $|\tilde a\tilde b|=|ab|$, $|\tilde a\tilde c|=|ac|$ and $|\tilde b\tilde c|=|bc|$.}
\end{defn}

In fact, in \cite{BGP}, it is always assumed that $X$ is locally complete. Our first main result shows that, in Definition 0.1, if $X$ is locally complete then condition (1) can be much weaker.

\vskip2mm

\noindent {\bf Theorem A.} \label{thm:A} {\it A locally complete intrinsic metric space $X$ is an Alexandrov space with curvature $\geqslant k$, if around any $x\in X$ there is a neighborhood $U_x$ such that the following condition is satisfied:
	\vskip1mm
	
	\noindent{\rm (A)}\ \ For any $p\neq q$ in $U_x$, there is a neighborhood $U_{q,p}$ of $q$, such that for any distinct $r,r_1,r_2\in U_{q,p}$,
	\begin{align}
		\widetilde{\angle}_kprr_1+\widetilde{\angle}_kprr_2+\widetilde{\angle}_kr_1rr_2\leqslant2\pi.
	\end{align} }
\hskip 4.5mm Note that $U_{q,p}$ can be sufficiently small, so the four points $p, r,r_1,r_2$ together looks very ``thin''. Moreover, in condition (A), if
(0.2) is replaced with $\widetilde{\angle}_kpqr_1+\widetilde{\angle}_kpqr_2+\widetilde{\angle}_kr_1qr_2\leqslant2\pi$, then it needs a
uniformity of $U_{q,p}$, i.e. $U_{q,p}$ contains an open ball $B(q,c|pq|)$ with $c$ being a positive constant not depending on $q$, but the `locally complete' is no longer needed (see Remark 2.3 below).

An intrinsic metric space $X$ is said to be {\it locally geodesic} if around  any $x\in X$ there is a neighborhood $U_x$ such that any two distinct points in $U_x$ can  be joined by a minimal geodesic (i.e. shortest path).
We will denote by $[pq]$ a minimal geodesic between $p$ and $q$ in $X$, and by $[pq]^\circ$ the interior part of $[pq]$. Via Alexandrov's lemma (Lemma 2.1 below), one can conclude that if $X$ in Definition 0.1 is in addition locally geodesic then condition (1) is equivalent to
(\cite{BGP}):

\vskip2mm

\noindent {\rm (2)}\ \  To any $p\in U_x$ and $[qr]\subset U_x$ with $p\not\in [qr]$, we associate $\tilde p\in\mathbb S^2_k$ and $[\tilde q\tilde r]\subset \mathbb S^2_k$ with $|\tilde p\tilde q|=|pq|$, $|\tilde p\tilde r|=|pr|$ and $|\tilde q\tilde r|=|qr|$. Then, for ALL $s\in [qr]^{\circ}$ and $\tilde s\in [\tilde q\tilde r]$ with $|qs|=|\tilde q\tilde s|$,
$$|ps|\geqslant |\tilde p\tilde s|.\eqno{(0.3)}$$
Moreover, in condition (2), (0.3) can be replaced with the version of (0.1) restricted to  $\{s,p,q,r\}$, i.e.
$$\widetilde{\angle}_kpsq+\widetilde{\angle}_kpsr\leqslant\pi. \eqno{(0.4)}$$
In condition (2), the requirement `$p\not\in [qr]$' is just to avoid `$p=s$' in (0.4).

Our second main result shows that condition (2) can also be much weaker.

\vskip2mm

\noindent {\bf Theorem B.} \label{thm:A} {\it
A locally geodesic intrinsic metric space $X$ is an Alexandrov space with curvature $\geqslant k$, if around any $x\in X$ there is a neighborhood $U_x$  such that one of the followings is satisfied:

\vskip1mm

\noindent{\rm (B1)}\ \ To any $p\in U_x$ and $[qr]\subset U_x$, we associate $\tilde p\in\mathbb S^2_k$ and $[\tilde q\tilde r]\subset \mathbb S^2_k$ with $|\tilde p\tilde q|=|pq|$, $|\tilde p\tilde r|=|pr|$ and $|\tilde q\tilde r|=|qr|$. Then, for $s\in [qr]$ and $\tilde s\in [\tilde q\tilde r]$ with $|qs|=|\tilde q\tilde s|$,
$$\limsup_{s\to q}\frac{|ps|-|\tilde p\tilde s|}{|qs|}\geqslant0.\eqno{(0.5)}$$

\vskip1mm

\noindent{\rm (B2)}\ \ For any small $\epsilon>0$, there is $\delta>0$ such that any $p\in U_x$ and $[qr]\subset U_x$  with $p\neq q$ satisfy
$$|ps|\leqslant|\bar p\bar s|+o(|qs|^2) \text{ with $o(|qs|^2)\leqslant\epsilon|qs|$ for all $s\in [qr]$ with $|qs|\leqslant\delta$},  \eqno{(0.6)}$$
where $\bar p$ and $\bar s$ belong to a triangle $\triangle\bar p\bar q\bar s\subset\mathbb S^2_k$ with $|\bar p\bar q|=|pq|$,
$|\bar q\bar s|=|qs|$ and $\angle\bar p\bar q\bar s=\limsup\limits_{t\to q,\ t\in [qr]}\tilde\angle_k pqt$.

\vskip1mm

\noindent{\rm (B3)}\ \ Given  $p\neq q$ in $U_x$, %there exists a neighborhood $U_{q,p}$ of $q$ such that
if $q$ lies in some $[r_1r_2]^\circ\subset U_x$ then, with respect to $t\triangleq\max\{|r_1q|, |qr_2|\}$,
$$\widetilde{\angle}_kpqr_1+\widetilde{\angle}_kpqr_2\leqslant\pi+o(t).\eqno{(0.7)}$$}
\hskip4.5mm For an equivalent version and a bit stronger version of condition (B2) and the relations between conditions (B1) and (B2), please refer to Section 1.

Recall that a Riemannian manifold $M$ is of sectional curvature $\geqslant k$ if and only if around any point in $M$ there is a neighborhood satisfying condition (2),
which is just the local version of the classical Toponogov's Theorem, and the fundamental tool of its proof is the second variation formula.
On the other hand, note that a bit stronger version of (0.5) is that $|ps|\geqslant |\tilde p\tilde s|+o(|qs|)$ for $s$ sufficiently close to $q$,
which and (0.6) are clearly related to the second variation formula of geodesic variation.
Thereby, if we define Alexandrov spaces with lower curvature bound using condition (B1) or (B2), then Theorem B can be viewed as a local version of Toponogov's Theorem on such spaces.

A significant property of an Alexandrov space $X$ with curvature $\geqslant k$ is that: {\it if $X$ is complete or geodesic \footnote{$X$ is said to be {\it geodesic} if
any two distinct points in $X$ can  be joined by a minimal geodesic.}, then {\rm (0.1)} holds for all $a,b,c,d\in X$}, which is called Globalization Theorem (\cite{BGP}, \cite{Petr2}).
There are several good tries to look for easier proofs for  the Globalization Theorem under the `complete and geodesic' condition (cf. [Pl], [Sh], [HSW]).
Inspired by the idea and using the key technique of the proofs of Theorems A and B, we can present new (and maybe more elementary) proofs for  Globalization Theorem only under `complete' and `geodesic'
respectively (note that `locally complete' and `locally geodesic' are respectively crucial to Theorems A and B).

\begin{remark}\label{rem0.2} {\rm In condition (B2), the infinitesimal $o(|qs|^2)$ (as $s\to q$) has to satisfy the requirement of the uniformity, i.e. $o(|qs|^2)\leqslant\epsilon|qs|$ when $|qs|\leqslant\delta$. As a counterexample, one can consider the space of three rays starting from a common point (which does not have lower curvature bound).}
\end{remark}

%\begin{remark}\label{rem0.3} {\rm In condition (B3), if any $U_{p,q}$ satisfies a much stronger property that, for any $[r_1r_2]\subset U_{p,q}$ and any $r_0\in [r_1r_2]^\circ$,
%$\widetilde{\angle}_kpr_0r_1+\widetilde{\angle}_kpr_0r_2\leqslant\pi$ (i.e., (0.2) restricted to $r_0,r_1,r_2$), then via Alexandrov's lemma it is easy to see that $U_x$ satisfies %condition (2) and thus is of curvature $\geqslant k$ (see Remark 4.7 below).}
%\end{remark}

\begin{remark}\label{rem0.4} {\rm As condition (A), condition (B1) has its corresponding `local' version (i.e.
it just need to consider all $[qr]$ in a small neighborhood $U_{q,p}$ around $q$, under which Theorem B is still true).
As for conditions (B2) and (B3), they are essentially just their `local' versions.}
\end{remark}

\begin{remark}\label{rem0.5} {\rm Since conditions (A) and (B1-3) are respectively weaker than conditions (1) and (2),
it will be possibly easier to check whether an intrinsic metric space satisfies the new conditions or not.
For instance, we can re-prove the Doubling Theorem by Perel$'$man in this way (see Section 5).}
\end{remark}

\begin{remark}\label{rem0.6} {\rm A locally geodesic intrinsic metric space $X$ is called an Alexandrov space with curvature $\leqslant k$ if
$$\text{ `$|ps|\geqslant |\tilde p\tilde s|$'\quad $\longrightarrow$\quad `$|ps|\leqslant |\tilde p\tilde s|$'\quad in condition (2)\ ([AKP]).}$$
Consequently and similarly, Theorem B has its `curvature $\leqslant k$' version, i.e., $X$ will be an Alexandrov space with curvature $\leqslant k$ if
we make the following changes in Theorem B:
$$\text{`$\limsup$' $\longrightarrow$ `$\liminf$',\quad `$\epsilon$' $\longrightarrow$ `$-\epsilon$' \text{ in (0.6)}, \quad `$\geqslant$' $\longleftrightarrow$ `$\leqslant$' except the `$\leqslant$' of `$|qs|\leqslant\delta$' in (0.6).}$$
In order to prove this, we just need to make the following changes in the corresponding arguments in Section 4:
$\text{`$\limsup$' $\longleftrightarrow$ `$\liminf$', `$+\epsilon$' $\longrightarrow$ `$-\epsilon$',  `$\geqslant$' $\longleftrightarrow$ `$\leqslant$', and `$\max$' $\longleftrightarrow$ `$\min$'.}$}
\end{remark}

\vskip2mm

In the rest of the paper, we will first show that condition (B1) is an ideal case of (B2) in Section 1. Then we will prove Theorem A in Section 2.
In Section 3, we will introduce a concept of `$f_\pm''(t)\leqslant B$ in the support sense', and then we prove Theorem B in Section 4.
As applications, we will supply new proofs for the Doubling Theorem in Section 5, and for the Globalization Theorem under the `complete' and `geodesic' conditions respectively in Section 6.

%%%%%%%%%%%%%%%%%%%%%%%%%%%%%%%%%%%%% section 1 relations %%%%%%%%%%%%%%%%%%%%%%%%%%%%%%%%%%%%%%%%%

\section{Relations between conditions (B1) and (B2)}

In the present paper, we will use the following model functions (cf. \cite{Pet}):
\begin{align*}
\sn_k(\rho)\triangleq
\begin{cases}\frac{1}{\sqrt{k}}\sin(\sqrt{k}\rho), & k>0 \\
\rho, & k=0 \\
\frac{1}{\sqrt{-k}}\sinh(\sqrt{-k}\rho), & k<0
\end{cases},\ \text{ct}_k(\rho)\triangleq\frac{\text{sn}_k'(\rho)}{\text{sn}_k(\rho)},\
f_k(\rho)\triangleq
\begin{cases}\frac{1}{k}\left(1-\cos(\sqrt{k}\rho)\right), & k>0 \\
\frac{1}{2}\rho^2, & k=0 \\
\frac{1}{-k}\left(\cosh(\sqrt{-k}\rho)-1\right), & k<0
\end{cases}.
\end{align*}

\subsection{An equivalent version of condition (B2)}

Condition (B2) can be formulated alternatively:

\noindent{\rm (B2)$'$}\ \ {\it For any small $\epsilon>0$, there is $\delta>0$ such that any $p\in U_x$ and $[qr]\subset U_x$  with $p\neq q$ satisfy
$$-\cos\tilde\angle_k pqs\leqslant -\cos\sphericalangle pqr+o(|qs|) \text{ with $o(|qs|)\leqslant\epsilon$ for all $s\in [qr]$ with $|qs|\leqslant\delta$}, \eqno{(1.1)}$$
where $\sphericalangle pqr\triangleq\limsup\limits_{t\to q,\ t\in [qr]}\tilde\angle_k pqt$.}

In fact, by the Law of Cosines on $\Bbb S_k^2$, for $s$ and $\bar s$ in (0.6) we have that
$$\begin{aligned}
|ps|&=|pq|-\cos\tilde\angle_k pqs \cdot |qs|+\frac12\text{ct}_k(|pq|)\sin^2\tilde\angle_k pqs \cdot |qs|^2+o_1(|qs|^2),\\
|\bar p\bar s|&=|pq|-\cos\sphericalangle pqr \cdot |qs|+\frac12\text{ct}_k(|pq|)\sin^2\sphericalangle pqr \cdot |qs|^2+o_2(|qs|^2),
\end{aligned}\eqno{(1.2)}
$$
where each $o_i(|qs|^2)$ is a higher order infinitesimal of $|qs|^2$ as $|qs|\to0$; or equivalently,
$$\begin{aligned}
f_k(|ps|)&=f_k(|pq|)-\text{sn}_k(|pq|)\cos\tilde\angle_k pqs \cdot |qs|+\frac12(1-kf_k(|pq|))\cdot |qs|^2+o_3(|qs|^2),\\
f_k(|\bar p\bar s|)&=f_k(|pq|)-\text{sn}_k(|pq|)\cos\sphericalangle pqr \cdot |qs|+\frac12(1-kf_k(|pq|))\cdot |qs|^2+o_4(|qs|^2).
\end{aligned}\eqno{(1.3)}
$$
Via (1.3), it is easy to see that (0.6) is equivalent to (1.1), i.e. condition (B2) is equivalent to (B2)$'$.

\subsection{Condition (B1) is an ideal case of condition (B2)}

Let $p,q,r,s$ and $\tilde p,\tilde q,\tilde r,\tilde s$ be the notations in condition (B1). Note that (0.5) holds automatically if $p=q$.
If $p\neq q$, by the Law of Cosines on $\Bbb S_k^2$ we have that
\begin{align*}
|ps|& =|pq|-\cos\tilde\angle_kpqs \cdot  |qs|+o_1(|qs|),\\
|\tilde p\tilde s|& =|pq|-\cos\tilde\angle_kpqr \cdot  |qs|+o_2(|qs|)
\end{align*}
for $s\in [qr]$ sufficiently close to $q$ and $\tilde s\in [\tilde q\tilde r]$ with $|qs|=|\tilde q\tilde s|$.
It then is easy to see that (0.5) in condition (B1) is equivalent to that $-\cos\tilde\angle_kpqr\leqslant-\liminf\limits_{s\to q}\cos\tilde\angle_kpqs$. Thereby, by replacing $[qr]$ with $[qs]$ for any $s\neq q$ in condition (B1), we have that
$$-\cos\tilde\angle_kpqs\leqslant-\cos\left(\limsup_{t\to q,\ t\in [qr]}\tilde\angle_k pqt\right),\ \text{ or equivalently},\ \tilde\angle_k pqs\leqslant \limsup_{t\to q,\ t\in [qr]}\tilde\angle_k pqt,\eqno{(1.4)}$$
which obviously fits (1.1). Namely, condition (B1) is essentially an ideal case of (B2)$'$, so of (B2).

\begin{remark}\label{rem1.1} {\rm Let $U_x$ satisfy condition (B1), and let $p\in U_x$ and $[r_1r_2]\subset U_x$  with $p\not\in [r_1r_2]$.
Since condition (B1) is a special case of (B2), by Lemma \ref{lem4.1} below we have that $\sphericalangle pqr_1+\sphericalangle pqr_2\leqslant\pi$
for any $q\in [r_1r_2]^\circ$,
where $\sphericalangle pqr_i\triangleq\limsup\limits_{t\to q,\ t\in [qr_i]}\tilde\angle_k pqt$. This plus (1.4) implies that
$$\tilde\angle_k pqr_1+\tilde\angle_k pqr_2\leqslant\pi$$
(without involving Lemma \ref{lem4.2}). I.e., $U_x$ satisfies (0.4) in condition (2).}
\end{remark}

%in Definition 0.1 by Alexandrov's lemma (Lemma 2.5 in \cite{BGP}, cf. \cite{AKP}).
%If it is under the local version of condition (A) (see Remark \ref{rem0.2}), one just only need to repeat Alexandrov's lemma finite times to see that $U_x$ satisfies condition (1) in Definition 0.1

\subsection{Other equivalent versions of condition (2)}

It is obvious that condition (2) has an equivalent version as  follows ([BGP]):

\vskip2mm

\noindent{\rm (3)}\ \ For any $p\in U_x$ and $[qr]\subset U_x$ with $p\neq q$, $\tilde\angle_k pqs$ with $s\in [qr]$ is decreasing with respect to $|qs|$.

\vskip2mm

\noindent Apparently, condition (3) implies that $\lim\limits_{t\to q,\ t\in [qr]}\tilde\angle_k pqt$ exists and $$\tilde\angle_k pqs\leqslant \lim_{t\to q,\ t\in [qr]}\tilde\angle_k pqt\ \text{ for ALL } s\in [qr]\setminus\{q\}.\eqno{(1.5)}$$

According to condition (3), for any $[pq]$ and $[qr]$ in an Alexandrov space
with curvature $\geqslant k$, we can define an angle between them at $q$ to be $\angle pqr\triangleq\lim\limits_{ x,\ y\to q}\tilde\angle_kxqy,$ where $x\in [pq]$ and $y\in[qr]$ (\cite{BGP}). It is clear that
$$\angle pqr\geqslant \widetilde{\angle}_k pqr \text{ for any $[pq],[qr]\subset U_x$}; \eqno{(1.6)}$$
or equivalently, given $[\bar p\bar q], [\bar q\bar r]\subset\mathbb S^2_k$ with $|\bar p\bar q|=|pq|$, $|\bar q\bar r|=|qr|$ and $\angle\bar p\bar q\bar r=\angle pqr$, we have that $|ps|\leqslant |\bar p\bar s|$
for all $s\in [qr]$ and $\bar s\in [\bar q\bar r]$ with $|\bar s\bar q|=|sq|$  (\cite{BGP}).
Moreover, if $q$ is an interior point of some $[rr']$ in addition, then $\angle pqr+\angle pqr'=\pi$.
Conversely, these properties together can be viewed as a sufficient condition for curvature $\geqslant k$. Namely,
condition (2)  can be replaced by the following conditions (cf. [BGP]):

\vskip1mm

\noindent{\rm (4-1)}\ \ For any $[pq]$ and $[qr]\subset U_x$, $\angle pqr$ can be defined so that
$\angle pqr+\angle pqr'\leqslant\pi$ if $q$ is in addition an interior point of some $[rr']$.

\noindent{\rm (4-2)}\ \ To any $[pq]$ and $[qr]\subset U_x$, we associate $[\bar p\bar q], [\bar q\bar r]\subset\mathbb S^2_k$ with $|\bar p\bar q|=|pq|$, $|\bar q\bar r|=|qr|$ and $\angle\bar p\bar q\bar r=\angle pqr$. Then $|ps|\leqslant |\bar p\bar s|$ for ALL $s\in [qr]$ and $\bar s\in [\bar q\bar r]$ with $|\bar s\bar q|=|sq|$.

\vskip1mm

\noindent Here, condition (4-1) is necessary (hint: one can define angles on the space of three rays starting from a common point so that (4-2) is satisfied, but (4-1) not).

\subsection{Weaker versions of conditions (3) and (4-2)}

Inspired by (1.1) and (1.5), we can present a bit stronger version of condition (B2)$'$ (so of (B2)) but a much weaker version of condition (3):

\vskip1mm

\noindent{\rm (B2)$''$}\ \ {\it For any small $\epsilon>0$, there is $\delta>0$ such that any $p\in U_x$ and $[qr]\subset U_x$ with $p\neq q$ satisfy
$$\tilde\angle_k pqs\leqslant \limsup\limits_{t\to q,\ t\in [qr]}\tilde\angle_k pqt +o(|qs|) \text{ with $o(|qs|)\leqslant\epsilon$, for all $s\in [qr]$ with $|qs|\leqslant\delta$}.\eqno{(1.7)}$$}
It is obvious that (1.7) implies (1.1), but not vice versa.

Inspired by (0.6), we would like to provide a weaker version of (4-2):

\vskip1mm

\noindent{\rm (C)}\ \ {\it To any $[pq]$ and $[qr]\subset U_x$, we associate $[\bar p\bar q], [\bar q\bar r]\subset\mathbb S^2_k$ with $|\bar p\bar q|=|pq|$, $|\bar q\bar r|=|qr|$ and $\angle\bar p\bar q\bar r=\angle pqr$. Then for $s\in [qr]$ sufficiently close to $q$ and $\bar s\in[\bar q\bar r]$ with $|\bar q\bar s|=|qs|$,
$$|ps|\leqslant|\bar p\bar s|+o(|qs|^2). \eqno{(1.8)}$$}
(Note that $\angle pqr$ has not been defined in (0.6).)

\begin{remark}\label{rem1.2} {\rm In proving Theorem B under condition (B2), we will first show that condition (4-1) is satisfied (see Lemma \ref{lem4.1}),
which plus (1.8) (partial information of (0.6)) guarantees Lemma \ref{lem4.2}, and then complete the proof almost immediately. In other words,
a locally geodesic intrinsic metric space is an Alexandrov space with curvature $\geqslant k$ if it satisfies conditions (4-1) and (C).}
\end{remark}

%It turns out that, among these conditions, condition (B) should be the weakest one to define an Alexandrov space with lower curvature bound (see Subsections 1.2 and 1.4).

%%%%%%%%%%%%%%%%%%%%%%%%%%%%%%%%%%%%% section 2 Proof of Theorem A%%%%%%%%%%%%%%%%%%%%%%%%%%%%%%%%%%%%%%%%%

\section{Proof of Theorem A}

The main goal of this section is to prove Theorem A. In the situation of Theorem A, there might be no minimal geodesic between $q$ and $r$,
so we cannot consider the comparison angle $\tilde\angle_kpqt$ as in condition (B2) and (1.5). However, the intrinsicness of the metric tells us that,
for any $d\in (0,|qr|)$, there is a sequence of points $\{s_i\}_{i=1}^\infty$  such that $|qs_i|\to d$ and $|s_ir|\to |qr|-d$ as $i\to\infty$.
So, we can define an `imaginary' comparison angle for $p,q,r$ and $d\in (0,|qr|)$ as follows:
$$\omega_k[_p^{qr}](d)\triangleq\inf\{\lim\limits_{i\to\infty}\widetilde{\angle}_kpqs_i\ |\ \{s_i\}_{i=1}^\infty\subset X,\ \mathcal{E}_{d}^{qr}(s_i)\to 0\ \text{ as } i\to\infty \},$$
where $$\mathcal{E}_d^{qr}(s_i)\triangleq\max\{||qs_i|-d|,\ ||s_ir|-(|qr|-d)|\}.$$
Note that $\{s_i\}_{i=1}^\infty$ might contain no converging subsequence even though $\mathcal{E}_{d}^{qr}(s_i)\to 0\ \text{ as } i\to\infty$. And
it is not hard to check that the `error' function $\mathcal{E}$ is additive with respect to $d$, i.e.
$$\text{if $\mathcal{E}_{d_1}^{qr}(s)<\varepsilon_1$ and $\mathcal{E}_{d_2}^{sr}(t)<\varepsilon_2$ with $d_1+d_2<|qr|$, then $\mathcal{E}_{d_1+d_2}^{qr}(t)<\varepsilon_1+\varepsilon_2$.}\eqno{(2.1)}$$

We now consider condition (A).  For convenience, we will call the $U_{q,p}$ in condition (A) a {\it good } neighborhood (associated to $p,q$),
and by the continuity of the distance we can set
$$	R_p(q)\triangleq\max\{R\ |\ B(q, R)\text{ is a good neighborhood}\}.$$
Note that, by Alexandrov's lemma ([BGP]), the goodness of $B(q,R_p(q))$ guarantees that
$$\text{$\omega_k[_p^{r_1r_2}](d)\geqslant \widetilde{\angle}_kpr_1r_2$ for any $r_1,r_2\in B(q,R_p(q))$ and $0<d<|r_1r_2|$};\eqno{(2.2)}$$
in particular, similar to the monotonicity of angles in condition (3), for any $r$ around $q$ we have that
$$\text{$\omega_k[_p^{qr}](d_1)\geqslant\omega_k[_p^{qr}](d_2)$, where $0<d_1<d_2<\min\{R_p(q),|qr|\}$.}\eqno{(2.3)}$$

\begin{lemma}[Alexandrov's lemma]\label{lem2.1} Let $\triangle pqr$, $\triangle pqs$, $\triangle abc \subset\Bbb S^2_k$
(where $\triangle pqr$ and $\triangle pqs$ are joined to each other in an exterior way along $[pq]$) such that
$|ab|=|pr|$, $|ac|=|ps|$, $|bc|=|qr|+|qs|$, and $|ab|+|ac|+|bc|<\frac{2\pi}{\sqrt k}$ if $k>0$. Then
$\angle pqr+\angle pqs\leqslant \pi$ (resp. $\geqslant\pi$) if and only if  $\angle prq\geqslant\angle abc$ and  $\angle psq\geqslant\angle acb$
(resp. $\angle prq\leqslant\angle abc$ and  $\angle psq\leqslant\angle acb$).
\end{lemma}

In fact, the goodness of $B(q,R_p(q))$ implies a stronger version of (2.2).

\begin{lemma}\label{lem2.2} Let $X$ be the space in Theorem A, and let $x\in X$ and $U_x$ satisfy condition (A). Then for distinct
$p,q,r\in U_x$ and $d\in(0,\frac{R_{p,r}(q)}{3}]$ with $R_{p,r}(q)\triangleq\min\{R_p(q),|qr|\}$, we have $\omega_k[_p^{qr}](d)\geqslant\widetilde{\angle}_kpqr$.
\end{lemma}

We will first prove Theorem A by assuming Lemma 2.2, and then verify the lemma.

\vskip2mm

\noindent{\it Proof of Theorem A}. Let $U_x$ be the neighborhood satisfying condition (A). Then for any distinct points $p,q,r,s\in U_x$, we just need to show that $\widetilde{\angle}_kpqr+\widetilde{\angle}_kpqs+\widetilde{\angle}_krqs\leqslant2\pi$ (see Definition 0.1).

Set $l\triangleq\min\{\frac{R_{p,s}(q)}{3},$ $\frac{R_{r,s}(q)}{3}\}$. For any small $\varepsilon>0$, by Lemma 2.2 there is $\bar s$ with $\mathcal{E}_{l}^{qs}(\bar s)$ sufficiently small such that
$$\widetilde{\angle}_kpqs\leqslant\omega_k[_p^{qs}](l)<\widetilde{\angle}_kpq\bar s+\varepsilon,\quad \widetilde{\angle}_krqs\leqslant\omega_k[_r^{qs}](l)<\widetilde{\angle}_krq\bar s+\varepsilon.$$
Similarly, for $l^\prime\triangleq\min\{\frac{R_{p,r}(q)}{3},\frac{R_{\bar s,r}(q)}{3}\}$, there is $\bar r$  with $\mathcal{E}_{l'}^{qr}(\bar r)$ sufficiently small such that
$$\widetilde{\angle}_k\bar sqr\leqslant\omega_k[_{\bar s}^{qr}](l')<\widetilde{\angle}_k\bar{s}q\bar r+\varepsilon, \quad \widetilde{\angle}_kpqr\leqslant\omega_k[_p^{qr}](l')<\widetilde{\angle}_kpq\bar r+\varepsilon.$$
Note that we can assume that both $\bar r$ and $\bar s$ lie in $B(q, R_p(q))$, so by (0.2) we have that
$$\widetilde{\angle}_kpq\bar r+\widetilde{\angle}_kpq\bar s+\widetilde{\angle}_k\bar{r}q\bar s\leqslant2\pi.$$
It then follows that
$$\widetilde{\angle}_kpqr+\widetilde{\angle}_kpqs+\widetilde{\angle}_krqs<\widetilde{\angle}_kpq\bar r+\widetilde{\angle}_kpq\bar s+\widetilde{\angle}_k\bar{r}q\bar s+4\varepsilon\leqslant2\pi+4\varepsilon.$$
Hence, we can conclude that $\widetilde{\angle}_kpqr+\widetilde{\angle}_kpqs+\widetilde{\angle}_krqs\leqslant2\pi$ due to the arbitrariness of $\varepsilon$. \hfill $\Box$

\vskip2mm

%Then we will prove Lemma 2.2. Note that the lemma is obviously true for trivial cases where $|pq|+|pr|=|qr|$, or $|pq|+|qr|=|pr|$, or $|qr|+|rp|=|qp|$.

%\vskip2mm

\renewcommand{\proofname}{\it Proof of Lemma 2.2}\begin{proof} By (2.3), it suffices to show that $\omega_k[_p^{qr}](\frac{R_{p,r}(q)}{3})\geqslant\widetilde{\angle}_kpqr$.
If it is not true, then it is clear that $|pq|+|pr|>|qr|$, and for $0<\delta\ll \min\{|pq|+|pr|-|qr|, |qr|,R_p(r)\}$ we claim that there is $\{s_i\}_{i=1}^\infty\subset U_x$ with $s_1=q$ such that
$$|qs_i|+|s_ir|-|qr|<2\delta\text{ (which implies $s_i\neq p$)}\eqno{(2.4)}$$
and
$$\text{(i)\ \ $\sum_{i=1}^\infty\frac{R_{p,r}(s_i)}{3}<|qr|-\delta$, \ \ (ii)\ \ $\mathcal{E}_{R_{p,r}(s_i)/3}^{s_ir}(s_{i+1})<\frac{\delta}{2^i}$}.$$
Note that (i) implies that $R_{p,r}(s_i)\to 0$ as $i\to\infty$, and (i) and (ii) together implies that
$$	\sum_{i=m}^{\infty}|s_is_{i+1}|<\sum_{i=m}^{\infty}\left(\frac{R_{p,r}(s_i)}{3}+\frac{\delta}{2^i}\right)\to 0\quad \text{as } m\to\infty. \eqno{(2.5)}$$
So, it follows that $\{s_i\}$ is a Cauchy sequence, and thus $s_i$ has to converge to a point $\bar{s}\in U_x\setminus\{r\}$ as $i\to\infty$
because we can assume that $U_x$ is complete by the local completeness of $X$. However, (2.4) implies $p\neq \bar s$, and thus $R_{p,r}(s_i)\stackrel{i\to\infty}{\longrightarrow} R_{p,r}(\bar{s})>0$, a contradiction.

We now need only to verify the claim above. In fact, each $s_i$ in the claim needs also to satisfy
$\omega_k[_p^{s_ir}](\frac{R_{p,r}(s_i)}{3})<\widetilde{\angle}_kps_ir$. Note that it is clear that $\frac{R_{p,r}(s_1)}{3}<|qr|$,
and `$\omega_k[_p^{s_1r}](\frac{R_{p,r}(s_1)}{3})<\widetilde{\angle}_kps_1r$' is the premise of the claim. It follows that there is
$\{s_{2}^j\}_{j=1}^\infty\subset U_x$ with $\mathcal{E}_{R_{p,r}(s_1)/3}^{s_1r}(s_{2}^j)\stackrel{j\to\infty}{\longrightarrow} 0$ such that $\lim\limits_{j\to\infty}\widetilde{\angle}_kps_1s_{2}^j<\widetilde{\angle}_kps_1r$.
Next, we show that $s_2^j$ for sufficiently large $j$ can be chosen as $s_2$.
Note that we can assume that  $\mathcal{E}_{R_{p,r}(s_1)/3}^{s_1r}(s_{2}^j)<\frac12\delta$ for large $j$, which implies $|qs_2^j|+|s_2^jr|-|qr|<\delta$.
And by Lemma 2.1, `$\lim\limits_{j\to\infty}\widetilde{\angle}_kps_1s_{2}^j<\widetilde{\angle}_kps_1r$' implies that
$$\lim_{j\to\infty}\widetilde{\angle}_kps_{2}^js_1+\lim_{j\to\infty}\widetilde{\angle}_kps_{2}^jr>\pi.\eqno{(2.6)}$$
On the other hand, note that $R_p(s_{2}^j)$ is almost bigger than $\frac23R_{p,r}(s_1)$ for large $j$. Then
due to the goodness of $B\left(s_2^j,R_p(s_{2}^j)\right)\ni s_1$, we have that
$$\widetilde{\angle}_kps_{2}^js_1+\omega_k[_p^{s_{2}^jr}](\frac{R_{p,r}(s_{2}^j)}{3})+\omega_k[_{s_1}^{s_{2}^jr}](\frac{R_{p,r}(s_{2}^j)}{3})\leqslant 2\pi\eqno{(2.7)}$$
(ref. (0.2)). Note that $\lim\limits_{j\to\infty}\omega_k[_{s_1}^{s_{2}^jr}](\frac{R_{p,r}(s_{2}^j)}{3})=\pi$, then it follows from (2.7) that
$$	\lim_{j\to\infty}\widetilde{\angle}_kps_{2}^js_1+\limsup_{j\to\infty}\omega_k[_p^{s_{2}^jr}](\frac{R_{p,r}(s_{2}^j)}{3})\leqslant \pi.\eqno{(2.8)}$$
By (2.6) and (2.8),  $\omega_k[_p^{s_{2}^{j}r}](\frac{R_{p,r}(s_{2}^{j})}{3})<\widetilde{\angle}_kps_{2}^{j}r$, which implies $|s_{2}^{j}r|\geqslant R_p(r)$ by (2.2).
Note that `$\mathcal{E}_{R_{p,r}(s_1)/3}^{s_1r}(s_{2}^{j})<\frac12\delta$' implies that
$\frac{R_{p,r}(s_1)}{3}+|s_{2}^{j}r|\leqslant \frac\delta2+|s_1s_{2}^{j}|+|s_{2}^{j}r|<|qr|+\frac32\delta$, and thus $$\frac{R_{p,r}(s_1)}{3}+\frac{R_{p,r}(s_{2}^{j})}3<|qr|+\frac32\delta-\frac{2|s_{2}^{j}r|}3<|qr|-\delta$$
(here the second `$<$' is due to $|s_{2}^{j}r|\geqslant R_p(r)\gg\delta$). We then can set $s_{2}\triangleq s_{2}^{j}$ for a sufficiently large $j$. Step by step, we can similarly locate $s_{i}$ for $i\geqslant3$ such that $\omega_k[_p^{s_{i}r}](\frac{R_{p,r}(s_{i})}{3})<\widetilde{\angle}_kps_{i}r$ (which implies $|s_ir|\geqslant R_p(r)$ by (2.2)), and  $\mathcal{E}_{R_{p,r}(s_{i-1})/3}^{s_{i-1}r}(s_i)<\frac1{2^{i-1}}\delta$ which plus (2.1) implies that $\mathcal{E}_{R_{p,r}(s_1)/3+\cdots+R_{p,r}(s_{i-1})/3}^{qr}(s_i)<\delta$
and so $|qs_i|+|s_ir|-|qr|<2\delta$, and that
$\frac{R_{p,r}(s_1)}{3}+\cdots+\frac{R_{p,r}(s_{i-1})}{3}+|s_{i}r|\leqslant \frac{\delta}2+\cdots+\frac{\delta}{2^{i-1}}+|s_1s_i|+|s_ir|<|qr|+3(\frac{\delta}2+\cdots+\frac{\delta}{2^{i-1}})$ and thus $$\frac{R_{p,r}(s_1)}{3}+\cdots+\frac{R_{p,r}(s_i)}3<|qr|+3\delta-\frac{2|s_{i}r|}3<|qr|-\delta.$$
That is, the claim above is verified.
\end{proof}

\begin{remark}\label{rem2.3} {\rm In condition (A), if (0.2) is replaced with $\widetilde{\angle}_kpqr_1+\widetilde{\angle}_kpqr_2+\widetilde{\angle}_kr_1qr_2\leqslant2\pi$,
then $U_{p,q}$ needs to contain a ball $B(q,c|pq|)$ with $c$ being a positive constant not depending on $q$, but the `locally complete' is no longer necessary.
In such a situation, all conclusions in the corresponding (2.2), (2.3) and Lemma 2.2 are still true, so is the corresponding Theorem A.
The proof of the corresponding Lemma 2.2 is almost a copy of that right above, where all $R_p(u)$ shall be replaced with $\frac12c|pu|$.
The difference is that, for some $i_0$, $B(s_{i_0},\frac12c|ps_{i_0}|)$ will contain $r$  (note that each $|ps_i|$ is bigger than a positive constant),
and thus a contradiction follows directly.}
\end{remark}

%%%%%%%%%%%%%%%%%%%%%%%%%%%%%%%%%%%%%%%%%%%%%%%% Section 3 %%%%%%%%%%%%%%%%%%%%%%%%%%%%%%%%%%%%%%%%%%%%%%%%%%%%%%%%%%%%%%

\section{ On\ \ $f_\pm''(t)+kf(t)\leqslant 0$ }

Recall that for a continuous function $f:(a,b)\to \mathbb{R}$, we say that $f^{\prime \prime}(t) \leqslant B$
(similarly for $f^{\prime \prime}(t) \geqslant B$) in the support sense at $t\in (a,b)$, if there is an $A \in \mathbb{R}$ such that
$$f(t+\tau) \leqslant\ f(t)+A \tau+\frac12B \tau^{2}+o\left(\tau^{2}\right) \text{ (cf. \cite{PP}, \cite{Petr1})}.\eqno{(3.1)}$$
And one can give a definition by a bit weaker conditions (refer to \S 10.5 in \cite{Na} and \S 9.3 in \cite{Pet})
\footnote{In \cite{Na} and \cite{Pet}, the definition is given respectively by: $\limsup\limits_{\tau\to 0}\frac{f(t+\tau)+f(t-\tau)-2f(t)}{\tau^2}\leqslant B$;
for any $\epsilon>0$,  there is a twice differentiable function $f_\epsilon$ defined on $(t-\delta,t+\delta)$ with some $\delta>0$  such that $f(t)=f_\epsilon(t),\ f(t+\tau)\leqslant\ f_\epsilon(t+\tau)\ \text{ for all } \tau\in(-\delta,\delta)$, and $f_\epsilon''(t)\leqslant\ B+\epsilon$.}. Note that `$f^{\prime \prime}(t) \leqslant B$' and `$f^{\prime \prime}(t)\geqslant B$' can hold simultaneously in the support sense, but $f$ does not have the second derivative at $t$ at all.

In this section, in order to prove Theorem B, we will introduce the right or left bounded second derivative in the support sense.

\vskip2mm

\noindent{\bf Definition 3.1.} Let $f: (a,b)\to \mathbb{R}$ be a continuous function.
At $t\in (a,b)$, we say that {\it $f_+^{\prime \prime}(t) \leqslant B$ in the support sense} if there is an $A \in \mathbb{R}$ such that
$$f(t+\tau) \leqslant f(t)+A \tau+\frac12B \tau^{2}+o\left(\tau^{2}\right) \text{ for } \tau>0,\eqno{(3.2)}$$
and $\liminf\limits_{\tau \to 0^-}\frac{f(t+\tau)-f(t)}{\tau}\geqslant A$; and `$f_-^{\prime \prime}(t) \leqslant B$' means that $f(t+\tau) \leqslant f(t)+A \tau+\frac12B \tau^{2}+o\left(\tau^{2}\right) \text{ for } \tau<0$ and  $\limsup\limits_{\tau \to 0^+}\frac{f(t+\tau)-f(t)}{\tau}\leqslant A$. Similarly, we can define $f_\pm^{\prime \prime}(t) \geqslant B$ in the support sense.

\vskip2mm

\noindent{\bf Remark 3.2.} It is easy to see that, in the support sense,  $f^{\prime \prime}(t) \leqslant B$ if and only if
$f_\pm^{\prime \prime}(t) \leqslant B$. And if $\limsup\limits_{\tau \to 0^+}\frac{f(t+\tau)-f(t)}{\tau}<\liminf\limits_{\tau \to 0^-}\frac{f(t+\tau)-f(t)}{\tau}$,
then $f_\pm^{\prime \prime}(t)=-\infty$ (i.e. $f_\pm^{\prime \prime}(t) \leqslant B$ for any $B\in \Bbb R$)  in the support sense.

\vskip2mm

\noindent{\bf Remark 3.3.}
It is obvious that $f^{\prime \prime}(t_0)\geqslant 0$  in the support sense if $f(t)$ achieves a local minimum at $t_0$. And note that if, for some $B_1$ and $B_2$,
$f_+^{\prime \prime}(t)\geqslant B_1$ and $f_+^{\prime \prime}(t)\leqslant B_2$ (or $f_-^{\prime \prime}(t)\geqslant B_1$ and $f_-^{\prime \prime}(t)\leqslant B_2$)
in the support sense, then $f'(t)$ exists and it has to hold that $B_1\leqslant B_2$.

\vskip2mm

In our proof of Theorem B, the key analysis tool is the following observation.

\vskip2mm

\noindent{\bf Lemma 3.4.} {\it  Let $f: [0,l]\to \mathbb{R}$ be a continuous function with $f(0)=f(l)=0$, and let $k$ be a constant real number such that $l<\frac{\pi}{\sqrt{k}}$ if $k>0$.  If
$$f_+''(t)+kf(t)\leqslant 0	\text{ or } f_-''(t)+kf(t)\leqslant 0\text{ in the support sense for all $t\in (0,l)$}, \eqno{\rm (3.3)}$$
then $f(t)\geqslant 0$ for all $t\in [0,l]$ \footnote{When $k>0$ and $l=\frac{\pi}{\sqrt{k}}$, one can consider $-\sin(\sqrt{k}t)$ as a counterexample.}.}

\vskip2mm

\noindent{\bf Remark 3.5.} In Lemma 3.4, the conclusion for $k=0$ is equivalent to: {\it a continuous function $f(t)|_{[a,b]}$ is concave if $f_+''(t)\leqslant 0$ or $f_-''(t)\leqslant 0$ in the support sense for all $t\in (a,b)$}. It is well known that $f(t)|_{[a,b]}$ is concave if $f''(t)\leqslant 0$ in the support sense for all $t\in (0,l)$ (\cite{Na});
so Lemma 3.4 for cases where $f''(t)+kf(t)\leqslant 0$ in the support sense for all $t\in (0,l)$ with $k\neq 0$ should be known to experts (cf. Lemma 1.8 in \cite{Petr1})
because the proof for cases where $k\neq0$ can be reduced to the case $k=0$.

\begin{proof}[Proof of  Lemma 3.4] Case 1: $k=0$. In this case, given $\epsilon>0$, we construct $$g(t)\triangleq f(t)-\frac{\epsilon}{2}t(t-l)$$
(which is the same as for `$f''(t)\leqslant 0$ in the support sense' in \cite{Na}).
It is clear that $g(t)$ is also continuous on $[0,l]$ with $g(0)=g(l)=0$,
and $g_+''(t)\leqslant-\epsilon$ or $g_-''(t)\leqslant-\epsilon$ in the support sense for all $t\in (0,l)$.
Observe that $g(t)\geqslant 0$ for all $t\in [0,l]$; otherwise,
$g(t)$ will achieve its (negative) minimum at some $t_0\in (0,l)$,
and thus $g''(t_0)\geqslant 0$ in the support sense which contradicts `$g_+''(t_0)\leqslant-\epsilon$ or $g_-''(t_0)\leqslant-\epsilon$' (see Remark 3.3).
It then follows that $f(t)\geqslant 0$ for all $t\in [0,l]$ when letting $\epsilon \to 0$.

\vskip1mm

Case 2: $k<0$. If $f(t)<0$ at some $t\in (0,l)$,
by the continuity there is $[a,b]\subseteq[0,l]$ such that $f(a)=f(b)=0$ and
$f(t)<0$ for all $t\in (a,b)$, and thus $f_+''(t)\leqslant -kf(t)<0$ or $f_-''(t)<0$ for all $t\in (a,b)$.
This is impossible by Case 1, i.e., it has to hold that $f(t)\geqslant 0$ for all $t\in [0,l]$.

\vskip1mm

Case 3: $k>0$. Due to the similarity, we just give a proof for $k=1$.
We argue by contradiction. Suppose that $f$ achieves its negative minimum at $t_0\in (0,l)$.
By the continuity of $f$, without loss of generality, we can assume that $f(t)\leqslant 0$ for all $t\in [0,l]$. Note that, for a positive number  $\lambda$, $\lambda f(t)$ still satisfies (3.3); so we can assume that $f(t_0)=-1$ by rescaling.

First of all, we will get a contradiction in the special case where $l\leqslant2\sqrt 2$. Note that we can assume that $l-t_0\leqslant t_0$ (or vice versa), so $l-t_0\leqslant\sqrt 2$.
For $t\in [t_0,l]$, we set
$$ h(t)\triangleq\frac{1}{(l-t_0)^2}(t-t_0)^2-1 \text{ and  } g(t)\triangleq f(t)-h(t).$$
Note that $g(t_0)=g(l)=0$. And since $h''(t)=\frac{2}{(l-t_0)^2}\geqslant 1$, we have that, in the support sense,
$$g_+''(t)= f_+''(t)-h''(t)\leqslant-f(t)-1\leqslant -f(t_0)-1=0\text{ or } g_-''(t)\leqslant 0.$$
Then by Case 1 and Remark 3.5, $g(t)$ is concave on $[t_0,l]$.
Moreover, note that `$f(t_0)=\min\{f(t)|\ t\in[0,l]\}$' and `$f_+''(t_0)\leqslant -f(t_0)=1$ or $f_-''(t_0)\leqslant 1$' together implies that
$f'(t_0)=0$ (see Remark 3.3), and thus $g'(t_0)=0$. This together with $g(t_0)=g(l)=0$ and the concavity of $g(t)$ implies that $g(t)\equiv 0$ on $[t_0,l]$.
It is impossible because $g_+''(t)\leqslant-f(t)-1<0$ or $g_-''(t)\leqslant-f(t)-1<0$ for $t$ near $l$.

We now can assume that $2\sqrt 2<l<\pi$, and still assume that $l-t_0\leqslant t_0$ (so $l-t_0<\frac\pi2$). Set
$$h(t)\triangleq-\sin(t-t_0+\frac\pi2) \text{ and  } g(t)\triangleq f(t)-h(t),\ t\in[t_0,l].$$
Note that $g(t_0)=0$ and $g(l)>0$, and $g_+''(t)+g(t)\leqslant 0$ or $g_-''(t)+g(t)\leqslant 0$ in the support sense on $[t_0,l]$.
Via the special case above (note that $l-t_0<\frac\pi2<2\sqrt 2$), we can conclude that $g(t)\geqslant 0$ on $[t_0,l]$, so $g_+''(t)\leqslant-g(t)\leqslant0$ or $g_-''(t) \leqslant0$.
Then $g(t)$ is concave on $[t_0,l]$ by Case 1 and Remark 3.5; moreover, we similarly have that $g'(t_0)=0$. Together with $g(t_0)=0$, these imply $g(t)\leqslant 0$ on $[t_0,l]$,
which contradicts $g(l)>0$.
\end{proof}

\noindent{\bf Remark 3.6.} Lemma 3.4 still holds if (3.3) is weakened to a piecewise version, namely,
$f_+''(t)+kf(t)\leqslant 0$ or $f_-''(t)+kf(t)\leqslant 0$  in the support sense for all but at most a finite number of $t\in (0,l)$, and
$\limsup\limits_{\tau \to 0^+}\frac{f(t+\tau)-f(t)}{\tau}\leqslant\liminf\limits_{\tau \to 0^-}\frac{f(t+\tau)-f(t)}{\tau}$ at each exceptional $t$.
This is obviously true in the special case where $k=0$ because, due to Remark 3.5, $f(t)$ is also concave in the situation here.
Then similar to Lemma 3.4, all other cases can be reduced to the special case.

\vskip2mm

\noindent{\bf Remark 3.7.} If the function $f(t)$ in Lemma 3.4 has the second derivative indeed, we would like to present an easier proof (cf. [Pet]).
Consider
$$g(t)\triangleq\frac{f(t)}{\sn_k(t)},\ t\in (0,l] \ (\text{note that } l<\frac{\pi}{\sqrt{k}} \text{ if } k>0),$$
and $g'(t)=\frac{f'(t)\sn_k(t)- f(t)\sn_k'(t)}{\sn_k^2(t)} \triangleq \frac{z(t)}{\sn_k^2(t)}$. Note that $\lim\limits_{t\to 0^+}z(t) =0$ and $z'(t)=(f''(t)+ kf(t))\sn_k(t)\leqslant 0$ by (3.3), which imply $z(t)\leqslant0$
and so $g'(t)\leqslant 0$ on $(0,l]$.  And thus, $g(t)\geqslant g(l)=0$, so does $f(t)$.

%%%%%%%%%%%%%%%%%%%%%%%%%%%%%%%%%%%%%%%%%%%%%%%%%%% Section 4%%%%%%%%%%%%%%%%%%%%%%%%%%%%%%%%%%%%%%%%%%%%%%%%%

\section{Proof of Theorem B}

We need only to prove Theorem B under condition (B2) and (B3) respectively because condition (B1) is an ideal case of (B2) (see Subsection 1.2).

\subsection{Under condition (B2)}

Recall that conditions (B2) and (B2)$'$ are equivalent to each other (see Subsection 1.1). First of all, we list two lemmas, and then prove Theorem B by assuming them.
\begin{lemma}\label{lem4.1}
Let $U_x$ be a neighborhood of $x\in X$ satisfying condition {\rm (B2)$'$}, and let $\gamma(t)|_{[0,\mu]}$ be an arc-length parameterized minimal geodesic in $U_x$.
Then for any $p\in U_x\setminus\gamma(t)|_{[0,\mu]}$ and $t_0\in(0,\mu)$,
$$\limsup\limits_{t\to t_0^+}\tilde\angle_k p\gamma(t_0)\gamma(t)+\limsup\limits_{t\to t_0^-}\tilde\angle_k p\gamma(t_0)\gamma(t)\leqslant\pi.\eqno{(4.1)}$$
\end{lemma}

\begin{lemma}\label{lem4.2}
Let $U_x$ be a neighborhood of $x\in X$ satisfying condition {\rm (B2)}, and let $\gamma(t)|_{[0,\mu]}$ be an arc-length parameterized minimal geodesic in $U_x$. Then for any $p\in U_x$, $g(t)\triangleq f_k(|p\gamma(t)|)$ satisfies
$$g''(t)+kg(t)\leqslant 1 \text{ in the support sense},\ \ \forall \ t\in(0,\mu).\eqno{(4.2)}$$
\end{lemma}

Refer to Section 1 for the function $f_k(\cdot)$. In proving Lemma \ref{lem4.2}, Lemma \ref{lem4.1} plays a crucial role.

\begin{proof} [Proof of Theorem B under condition (B2)]\

Let $U_x$ be a neighborhood of $x\in X$ satisfying condition (B2), and let $\gamma(t)|_{[0,\mu]}$ be an arc-length parameterized minimal geodesic in $U_x$. To $\gamma(t)|_{[0,\mu]}$ and any $p\in U_x$, we associate an arc-length parameterized minimal geodesic $\tilde\gamma(t)|_{[0,\mu]}$ and a point $\tilde p$ in $\mathbb S^2_k$ with
$|\tilde p\tilde\gamma(0)|=|p\gamma(0)|$ and $|\tilde p\tilde\gamma(\mu)|=|p\gamma(\mu)|$.
According to condition (2), it suffices to show that
$$|p\gamma(t)|\geqslant |\tilde p\tilde\gamma(t)|\ \ \forall\  t\in[0,\mu].\eqno{(4.3)}$$
We let $g(t)\triangleq f_k(|p\gamma(t)|)$ and $\tilde g(t)\triangleq f_k(|\tilde p\tilde\gamma(t)|)$.
Note that $g(0)=\tilde g(0)$ and $g(\mu)=\tilde g(\mu)$. Moreover, it is well known that (\cite{PP}, \cite{Pet})
$$\tilde g''(t)+k\tilde g(t)=1\ \ \forall \ t\in(0,\mu),$$
then by Lemma \ref{lem4.2} it is clear that, in the support sense,
$$\left(g(t)-\tilde g(t)\right)''+k\left(g(t)-\tilde g(t)\right)\leqslant 0\ \ \forall \ t\in(0,\mu).\eqno{(4.4)}$$
And thus, we can apply Lemma 3.4 to conclude that $g(t)-\tilde g(t)\geqslant 0$  for all $t\in[0,\mu]$, i.e. $f_k(|p\gamma(t)|)\geqslant f_k(|\tilde p\tilde\gamma(t)|)$, which is equivalent to (4.3). (Here, it needs that $\mu<\frac{\pi}{\sqrt{k}}$ if $k>0$, while this can be guaranteed indeed because $U_x$ can be sufficiently small.)
\end{proof}

\begin{remark}\label{rem2.3} {\rm \cite{PP} has mentioned (without proof) that: {\it a geodesic intrinsic metric space $X$ is an Alexandrov space with curvature $\geqslant k$ if and only if {\rm (4.2)} is satisfied for any $p\in X$ and any minimal geodesic $\gamma(t)\subset X$}. This coincides with the idea of proving Toponogov's Theorem on Riemannian manifolds (cf. \cite{Pet}). So, the argument in the proof right above should be almost obvious to experts. }
\end{remark}

Next, we verify Lemma \ref{lem4.1}, the key lemma.

\begin{proof} [Proof of Lemma \ref{lem4.1}] \

As in (B2)$'$, we let $\sphericalangle p\gamma(t_0)\gamma(\mu)\triangleq\limsup\limits_{t\to t_0^+}\tilde\angle_k p\gamma(t_0)\gamma(t)$ and $\sphericalangle p\gamma(t_0)\gamma(0)\triangleq\limsup\limits_{t\to t_0^-}\tilde\angle_k p\gamma(t_0)\gamma(t)$.

\noindent{\bf Claim}: $\sphericalangle p\gamma(t_0)\gamma(\mu)\leqslant\liminf\limits_{t\to t_0^-}\sphericalangle p\gamma(t)\gamma(\mu)$.

According to (1.1) in condition (B2)$'$, for any small $\epsilon>0$ there is $\delta>0$ such that
$$-\cos\tilde{\angle }_k p\gamma(t)\gamma(t+\tau)\leqslant-\cos\sphericalangle p\gamma(t)\gamma(\mu)+\epsilon$$
for all $t\in (0,\mu)$ and any $\tau\in (0,\delta)$ with $t+\tau\leqslant \mu$.
And, by the continuity of $|p\gamma(t)|$, it is clear that
$$\lim_{t\to t_0^-}\tilde{\angle }_k p\gamma(t)\gamma(t+\tau)=\tilde{\angle }_k p\gamma(t_0)\gamma(t_0+\tau).$$
It then follows that
$$-\cos\tilde{\angle }_k p\gamma(t_0)\gamma(t_0+\tau)\leqslant-\cos\left(\liminf\limits_{t\to t_0^-}\sphericalangle p\gamma(t)\gamma(\mu)\right)+\epsilon.$$
Note that the claim follows as long as we let $\tau\to 0$ first and then let $\epsilon\to 0$.

We now consider the distance function $f(t)\triangleq |p\gamma(t)|$, and for any $t_0\in (0, \mu)$ we set
\begin{align*}
f_{-,\min}^{\prime}(t_0)\triangleq\liminf_{t\to t_0^-}\frac{f(t)-f(t_0)}{t-t_0} \quad\text{and} \quad f_{+,\max}^{\prime}(t_0)\triangleq\limsup_{t\to t_0^+}\frac{f(t)-f(t_0)}{t-t_0}.
\end{align*}
Observe that
$$f_{-,\min}^{\prime}(t_0)=\cos\sphericalangle p\gamma(t_0)\gamma(0)\ \text{ and }\ f_{+,\max}^{\prime}(t_0)=-\cos\sphericalangle p\gamma(t_0)\gamma(\mu).$$
In fact, as $t\to t_0^-$, it is easy to see that
$$f(t)-f(t_0)=-\cos\tilde{\angle}_k p\gamma(t_0)\gamma(t)\cdot (t_0-t)+o(t_0-t)$$
(by the Law of Cosines on $\Bbb S_k^2$). And thus
$$f_{-,\min}^{\prime}(t_0)=\liminf_{t\to t_0^-}\left(\cos\tilde{\angle}_k p\gamma(t_0)\gamma(t)\right)=\cos\left(\limsup_{t\to t_0^-}\tilde{\angle}_k p\gamma(t_0)\gamma(t)\right)=\cos\sphericalangle p\gamma(t_0)\gamma(0).$$
Similarly, we have that $f_{+,\max}^{\prime}(t_0)=-\cos\sphericalangle p\gamma(t_0)\gamma(\mu)$.

Another important observation is that
$\liminf\limits_{t\to t_0^-}f_{+,\max}^{\prime}(t)\leqslant f_{-,\min}^{\prime}(t_0)$ (by the Sublemma below).
As a result, we have that
$$\cos\sphericalangle p\gamma(t_0)\gamma(0)\geqslant \liminf\limits_{t\to t_0^-}\left(-\cos\sphericalangle p\gamma(t)\gamma(\mu)\right)\geqslant-\cos\left(\liminf_{t\to t_0^-}\sphericalangle p\gamma(t)\gamma(\mu)\right),$$
which together with the claim above implies that
$$\cos\sphericalangle p\gamma(t_0)\gamma(0)\geqslant -\cos\sphericalangle p\gamma(t_0)\gamma(\mu),\ \ \text{or equivalently,} \ \ \sphericalangle p\gamma(t_0)\gamma(0)+\sphericalangle p\gamma(t_0)\gamma(\mu)\leqslant\pi.$$
This completes the proof.
\end{proof}

\noindent{\bf Sublemma.} {\it
Let $f: (a,b)\to \mathbb{R}$ be a continuous function. Then for any $t_0\in (a,b)$,
$$\liminf_{t\to t_0^-}f_{+,\max}^{\prime}(t)\leqslant f_{-,\min}^{\prime}(t_0).$$}
\noindent{\it Proof}. (This might be a known result in real analysis.)
We first note that, for any $[c,d]\subset(a,b)$, there exists $x\in [c,d)$ such that
$$f_{+,\max}^{\prime}(x)\triangleq \limsup_{t\to x^+}\frac{f(t)-f(x)}{t-x}
\leqslant \frac{f(d)-f(c)}{d-c}$$
(which is similar to the Differential Mean Value Theorem, cf. \cite{MV}). Meantime, by definition,
$$f_{-,\min}^{\prime}(t_0)=\liminf_{t\to t_0^-}\frac{f(t_0)-f(t)}{t_0-t}.$$
Therefore, there exists $x_t\in [t,t_0)$ (here $a<t<t_0$) such that
$$\hskip3.5cm\liminf_{t\to t_0^-} f_{+,\max}^{\prime}(t)\leqslant\liminf_{t\to t_0^-} f_{+,\max}^{\prime}(x_t)\leqslant f_{-,\min}^{\prime}(t_0).\hskip3.5cm\Box$$

\begin{remark}\label{rem2.4} {\rm Note that in proving Lemma \ref{lem4.1}, we use only the uniformity in (1.1),
i.e. ``$-\cos\tilde\angle_k pqs\leqslant -\cos\sphericalangle pqr+\epsilon$ for all $s\in [qr]$ with $|qs|\leqslant\delta$'', without involving $o(|qs|)$.}
\end{remark}

Now, in order to complete the proof of Theorem B under (B2), we just need to verify  Lemma \ref{lem4.2}.

\begin{proof} [Proof of Lemma \ref{lem4.2}]

Let $U_x$ satisfy condition (B2) (or (B2)$'$), and let $p\in U_x$ and $[qr]\subset U_x$ with $p\neq q$, and let $s, \bar p, \bar s$ be the notations in (0.6). By (0.6), it holds that $|ps|\leqslant |\bar p\bar s|+o(|qs|^2)$, which is equivalent to $f_k(|ps|)\leqslant f_k(|\bar p\bar s|)+o(|qs|^2)$. Then by the second equality of (1.3), we have that
$$f_k(|ps|)\leqslant f_k(|pq|)-\text{sn}_k(|pq|)\cos\sphericalangle pqr \cdot |qs|+\frac12(1-kf_k(|pq|))\cdot |qs|^2+o(|qs|^2),\eqno{(4.5)}$$
where $\sphericalangle pqr\triangleq\limsup\limits_{t\to q,\ t\in [qr]}\tilde\angle_k pqt$.
We now consider the $\gamma(t)|_{[0,\mu]}\subset U_x$ in the lemma. Note that we just need to consider the case where $p\not\in\gamma(t)|_{[0,\mu]}$.
Then for any $t_0\in (0,\mu)$, it holds that $\sphericalangle p\gamma(t_0)\gamma(0)+\sphericalangle p\gamma(t_0)\gamma(\mu)\leqslant\pi$ by Lemma 4.1.
And thus, via (4.5) it is easy to see that
$$f_k(|p\gamma(t_0+\tau)|)\leqslant f_k(|p\gamma(t_0)|)-\text{sn}_k(|p\gamma(t_0)|)\cos\sphericalangle p\gamma(t_0)\gamma(\mu) \cdot \tau+\frac12(1-kf_k(|p\gamma(t_0)|))\cdot \tau^2+o(\tau^2)$$
where $|\tau|$ is sufficiently small. Then, it is clear that $g''(t_0)\leqslant 1 -kg(t_0)$ in the support sense, where $g(t)=f_k(|p\gamma(t)|)$; i.e., the proof of the lemma is completed.
\end{proof}

\begin{remark}\label{rem4.5} {\rm We would like to provide a bit weaker version of condition (B2):
{\it For any $p\in U_x$ and $[r_1r_2]\subset U_x$ with $p\not\in[r_1r_2]$, {\rm(0.6)} holds with respect to $[qr_i]$ for all but at most a finite number of $q\in [r_1r_2]^\circ$, and {\rm(0.6)}
is replaced by ``$|ps|\leqslant|\bar p\bar s|+\epsilon|qs|$ for all $s\in [qr_i]$ with $|qs|\leqslant\delta$'' at each exceptional $q$ (cf. Remark 4.4).}
Inspired by Remark \ref{rem1.2}, one of its alternative version is: {\it {\rm(4-1)} holds on $U_x$;
and for any $p\in U_x$ and $[r_1r_2]\subset U_x$ with $p\not\in[r_1r_2]$, {\rm (1.8)} holds with respect to $[qr_i]$ for all but at most a finite number of $q\in [r_1r_2]^\circ$, and at each exceptional $q$ it holds that $|ps|\leqslant |pq|-\cos\angle pqr_i\cdot |qs|+o(|qs|)$ for $s\in [qr_i]$ sufficiently close to $q$}. %This guarantees that  $f'_{+,\max}(t_0)\leqslant f'_{-,\min}(t_0)$ if we set $f(t)\triangleq|p\gamma(t)|$, where $\gamma(t)|_{[0,|r_1r_2|]}$ denotes $[r_1r_2]$ with $|r_1\gamma(t)|=t$, and
As in proving Theorem B under condition (B2),
we can apply Remark 3.6 instead of Lemma 3.4 to see that $X$ is of curvature $\geqslant k$ under the weaker versions.}
\end{remark}

\subsection{Under condition (B3)}

For the same reason as in Subsection 4.1 (cf. (4.4)), we just need to establish the following lemma in order to give a proof for Theorem B under condition (B3).

\begin{lemma}\label{lem4.6}
Let $U_x$ be the neighborhood of $x\in X$ satisfying condition {\rm (B3)}, and let $\gamma(t)|_{[0,\mu]}$ be an arc-length parameterized minimal geodesic in $U_x$. Then for any $p\in U_x$, $g(t)\triangleq f_k(|p\gamma(t)|)$ satisfies
$$g_+''(t)+kg(t)\leqslant 1 \text{ or } g_-''(t)+kg(t)\leqslant 1 \text{ in the support sense},\ \ \forall \ t\in(0,\mu).\eqno{(4.6)}$$
\end{lemma}

Here, we can not derive that $g''(t)+kg(t)\leqslant 1$ in the support sense as in (4.2), although condition (B3) implies immediately
the inequality in (4.1), i.e. for any $t_0\in(0,\mu)$ in Lemma 4.6
$$\limsup\limits_{t\to t_0^+}\tilde\angle_k p\gamma(t_0)\gamma(t)+\limsup\limits_{t\to t_0^-}\tilde\angle_k p\gamma(t_0)\gamma(t)\leqslant\pi.\eqno{(4.7)}$$

\begin{proof} [Proof of Lemma \ref{lem4.6}]

We just need to consider the case where $p\not\in\gamma(t)|_{[0,\mu]}$. For $t_0\in (0,\mu)$,
via (1.3) we can conclude that
\begin{align*}&\limsup\limits_{t\to t_0^+}\frac{g(t)-g(t_0)}{t-t_0}=-\text{sn}_k(|p\gamma(t_0)|)\cos\sphericalangle p\gamma(t_0)\gamma(\mu) \text{ and } \\
& \liminf\limits_{t\to t_0^-}\frac{g(t)-g(t_0)}{t-t_0}=\text{sn}_k(|p\gamma(t_0)|)\cos\sphericalangle p\gamma(t_0)\gamma(0),
\end{align*}
where $\sphericalangle p\gamma(t_0)\gamma(0)=\limsup\limits_{t\to t_0^-}\tilde\angle_k p\gamma(t_0)\gamma(t)$ and
$\sphericalangle p\gamma(t_0)\gamma(\mu)=\limsup\limits_{t\to t_0^+}\tilde\angle_k p\gamma(t_0)\gamma(t)$. So, by (4.7)
$$\limsup\limits_{t\to t_0^+}\frac{g(t)-g(t_0)}{t-t_0}\leqslant\liminf\limits_{t\to t_0^-}\frac{g(t)-g(t_0)}{t-t_0}.$$
This plus Remark 3.2 guarantees that (4.6) holds at $t_0$ if we can show that at least one of the following holds:
\begin{align*} & g(t_0+\tau)\leqslant g(t_0)-\text{sn}_k(|p\gamma(t_0)|)\cos\sphericalangle p\gamma(t_0)\gamma(\mu) \cdot \tau+\frac12(1-kg(t_0))\cdot \tau^2+o(\tau^2) \text{ for }  \tau>0;\\
\ \ \ \ & g(t_0+\tau)\leqslant g(t_0)+\text{sn}_k(|p\gamma(t_0)|)\cos\sphericalangle p\gamma(t_0)\gamma(0) \cdot \tau+\frac12(1-kg(t_0))\cdot \tau^2+o(\tau^2)\text{ for }  \tau<0;\ \ \ \ \  (4.8)\\
&-\cos\sphericalangle p\gamma(t_0)\gamma(\mu)<\cos\sphericalangle p\gamma(t_0)\gamma(0).
\end{align*}
If any of the three inequalities is not true, then $-\cos\sphericalangle p\gamma(t_0)\gamma(\mu)=\cos\sphericalangle p\gamma(t_0)\gamma(0)$ (i.e. $\sphericalangle p\gamma(t_0)\gamma(0)+\sphericalangle p\gamma(t_0)\gamma(\mu)=\pi$), and by taking into account the first equality in (1.3) there is $c>0$, $\tau_{+,i}>0$ and $\tau_{-,i}<0$ with $\tau_{+,i},\tau_{-,i}\to0$ as $i\to\infty$ such that
$$\widetilde{\angle}_kp\gamma(t_0)\gamma(t_0+\tau_{+,i})>\sphericalangle p\gamma(t_0)\gamma(\mu)+c\tau_{+,i}\text{ and }
\widetilde{\angle}_kp\gamma(t_0)\gamma(t_0+\tau_{-,i})>\sphericalangle p\gamma(t_0)\gamma(0)-c\tau_{-,i}.$$
It then follows that $\widetilde{\angle}_kp\gamma(t_0)\gamma(t_0+\tau_{+,i})+\widetilde{\angle}_kp\gamma(t_0)\gamma(t_0+\tau_{-,i})>\pi+c\max\{\tau_{+,i},-\tau_{-,i}\}$,
which contradicts (0.7) in condition (B3).
\end{proof}

\noindent{\bf Remark 4.7.} If (0.7) in condition (B3) is strengthened to `$\widetilde{\angle}_kpqr_1+\widetilde{\angle}_kpqr_2\leqslant\pi$' with $r_i$ sufficiently close to $q$, then
in the proof right above it is clear that the former two inequalities of (4.8) hold and
$\sphericalangle p\gamma(t_0)\gamma(0)+\sphericalangle p\gamma(t_0)\gamma(\mu)\leqslant\pi$,
where $\sphericalangle p\gamma(t_0)\gamma(0)=\sup\limits_{t_0-\tau_{t_0}<t<t_0}\widetilde\angle_kp\gamma(t_0)\gamma(t)$ for sufficiently small $\tau_{t_0}>0$
and similar for $\sphericalangle p\gamma(t_0)\gamma(\mu)$. It then follows that  $g''(t_0)+kg(t_0)\leqslant 1$ in the support sense.

%%%%%%%%%%%%%%%%%%%%%%%%%%%%%%%%%%%%% section 5  New proof  %%%%%%%%%%%%%%%%%%%%%%%%%%%%%%%%%%%%%%%%%

\section{New proof of the Doubling Theorem}

As mentioned in Remark 0.4,  it will be possibly easier to check whether an intrinsic metric space satisfies the new conditions (A-C) or not.
For instance, one can check that the doubling space in the Doubling Theorem (\cite{Pe})  satisfies the weaker version of condition (B2) in Remark 4.5. In other words, we can give a new and more detailed proof for the theorem.

\vskip2mm

\noindent{\bf Theorem 5.1} (Doubling Theorem). {\it Let $X$ be a complete $n$-dimensional Alexandrov space with curvature $\geqslant k$ and
with nonempty boundary. Then its doubling $D(X)$ with canonical metric is  a complete $n$-dimensional Alexandrov space with curvature $\geqslant k$ and with empty boundary.}

\vskip2mm

We can view $D(X)$ as $X\cup_{\partial X}\bar X$, where $\bar X$ is a copy of $X$. Between any $x$ and $y$ in $D(X)$, the distance $|xy|$ is just the distance between them in $X$ (or $\bar X$) if both $x$ and $y$ lie in $X$ (or $\bar X$), otherwise, $|xy|\triangleq\min\limits_{z\in\partial X}\{|xz|+|zy|\}$ (\cite{Pe}). Then it is easy to see the following property of $D(X)$.

\vskip2mm

\noindent{\bf Fact 5.2} (cf. \cite{Pe}). {\it For any $[xy]\subset D(X)$, $[xy]^\circ$
passes through $\partial X$ at most once or $[xy]\subset\partial X$ ($\subset D(X)$), and $[xy]^\circ$ passes through $\partial X$ once if and only if $x\in X^\circ$ and $y\in \bar X^\circ$ or vice versa.}

\vskip2mm

In the proof of Theorem 5.1, we will use the following fundamental property of $X$.
\vskip2mm

\noindent{\bf Fact 5.3} (cf. \cite{BGP}). {\it For $q\in \partial X$ and small $\varepsilon>0$,   $\max\limits_{\zeta\in\partial\Sigma_qX}\min\{|\zeta\uparrow_q^u|:\ u\in \partial X,\ |qu|=\varepsilon\}\to0$ as $\varepsilon\to0$ \footnote{An essential reason for this is that
$(\frac{1}{\varepsilon}X, q)\to T_q$ with $(\frac{1}{\varepsilon}\partial X, q)\to \partial T_q$ as $\varepsilon\to0$, where $T_q$ is the tangent cone of $X$ at $q$ (\cite{BGP}).
In fact, $\max\limits_{\zeta\in\Sigma_qX}\min\{|\zeta\uparrow_q^u|:\ u\in X,\ |qu|=\varepsilon\}\to 0$ as $\varepsilon\to0$ whether $q\in\partial X$ or not. }, where $\Sigma_qX$ is the space of directions of $X$ at $q$
and $\uparrow_q^u\in\Sigma_qX$ is the direction of some $[qu]$. And for any $[qx]$ with $x\in X^\circ$, $\uparrow_q^x\in\left(\Sigma_qX\right)^\circ$.}

\vskip2mm

\begin{proof}[Proof of Theorem 5.1]

Note that we can apply induction on the dimension $n$ starting with the trivial case where $n=1$
(it is our convention that $X$ is an interval of length $\leqslant\pi$ if $n=1$, cf. \cite{BGP}).

We first show that $D(X)$ satisfies condition (4-1) globally; i.e., for any $[pq]$ and $[qr]\subset D(X)$, $\angle pqr$ can be defined so that $\angle pqr+\angle pqr'\leqslant \pi$ if $q$ is an interior point of some $[rr']$.
By induction, for any $x\in \partial X$, $D(\Sigma_xX)=\Sigma_xX\cup_{\partial\Sigma_xX}\Sigma_x\bar X$
is  a complete $(n-1)$-dimensional Alexandrov space with curvature $\geqslant 1$ and with empty boundary.
Then based on Fact 5.2, we can define $\angle pqr$ to be the distance between the directions of $[pq]$ and $[qr]$ at $q$ in
$D(\Sigma_qX)$ or $\Sigma_qX$ according to $q\in\partial X$ or not respectively. And if $q\in [rr']^\circ$, there are only the following possible cases.

$\bullet$ $q\in X^\circ$ or $\bar X^\circ$, or $[rr']\subset\partial X$: It is clear that $\angle pqr+\angle pqr'=\pi$ by (4-1) on $X$ (or $\bar X$).

$\bullet$ $q\in \partial X$, and  $r\in X^\circ$ and $r'\in \bar X^\circ$ (or vice versa): In this case, any $[qr]$ lies in $X$ and $[qr']$ lies in $\bar X$,
and $\uparrow_q^r\in\left(\Sigma_qX\right)^\circ$ and $\uparrow_q^{r'}\in\left(\Sigma_q\bar X\right)^\circ$ (see Facts 5.2 and 5.3). By induction on $D(\Sigma_qX)$ and the first variation formula on $X$ and $\bar X$, it has to hold that $|\uparrow_q^r\uparrow_q^{r'}|=\pi$ in $D(\Sigma_qX)$. As a result, $\angle pqr+\angle pqr'=\pi$, and there is a unique minimal geodesic between $q$ and $r$ (or $r'$).

\vskip1mm

Next, we let $[pq]$ and $[qr]$ be two minimal geodesics in $D(X)$ with $p\in\bar X$, and
$|pq|+|qr|<\frac{\pi}{\sqrt k}$ if $k>0$. Let $\triangle\bar p\bar q\bar r$ be a triangle in $\mathbb S^2_k$ with $|\bar p\bar q|=|pq|$, $|\bar q\bar r|=|qr|$ and $\angle\bar p\bar q\bar r=\angle pqr$, and let $s\in [qr]$ and $\bar s\in [\bar q\bar r]$ satisfy $|\bar s\bar q|=|sq|$.

\noindent\hskip1.5mm{\bf Claim 1}. {\it If $q\in \partial X$, then for $s\in [qr]$ sufficiently close to $q$ we have that $$|ps|\leqslant |\bar p\bar s|+o(|qs|).\eqno{(5.1)}$$}
{\bf Claim 2}. {\it If $q\in X^\circ$, then for $s\in [qr]$ sufficiently close to $q$ we have that $$|ps|\leqslant |\bar p\bar s|+o(|qs|^2).\eqno{(5.2)}$$}
\hskip4.5mm Note that (5.2) holds naturally by condition (4-2) on $\bar X$ if $q\in \bar X^\circ$. Then Claims 1 and 2 together with Fact 5.2 and that $D(X)$ satisfies (4-1) enable us to apply Remark \ref{rem4.5} to conclude that $D(X)$ is of curvature $\geqslant k$, and then it follows that  $\Sigma_xD(X)=D(\Sigma_xX)$ or $\Sigma_xX$ according to $x\in\partial X$ or not respectively, which implies that $D(X)$ has no boundary point. Namely, the Doubling Theorem follows (note that the completeness of $D(X)$ is obvious). So, it remains to verify the two claims.

\vskip1mm

We verify Claim 1 at first. By (4-2) on $X$ or $\bar X$, (5.2) is true, so is (5.1), if $p\in\partial X$ or $[qr]\subset\bar X$. I.e., we can assume that $p\in \bar X^\circ$ (and so $\uparrow_q^p\in (\Sigma_q\bar X)^\circ$) and $\uparrow_q^r\in (\Sigma_q X)^\circ$.
Then there is $\xi\in\partial\Sigma_qX$ such that $|\uparrow_q^p\xi|+|\xi\uparrow_q^r|=|\uparrow_q^p\uparrow_q^r|\ (=\angle pqr)$ in $D(\Sigma_qX)$. Here, the difficulty is that there might be no $[qu]\subset\partial X$ with $\uparrow_q^u=\xi$. Nevertheless, according to Fact 5.3, for sufficiently small $\epsilon>0$ there is $\varepsilon_\epsilon$ such that for any positive $t\leqslant\varepsilon_\epsilon$ there is $u\in\partial X$ with $|qu|=t$ and $|\uparrow_q^u\xi|\leqslant\epsilon$ in $\Sigma_qX$. Note that we can assume that $\angle pqr<\pi$ because (5.1) is obviously true if $\angle pqr=\pi$.
Then for $s\in[qr]$ sufficiently close to $q$, there is $\bar u\in [\bar s\bar p]$ such that $\angle\bar p\bar q\bar u=|\uparrow_q^p\xi|$ (and then $\angle\bar s\bar q\bar u=|\uparrow_q^s\xi|=|\uparrow_q^r\xi|$), and thus there is $u\in\partial X$ such that $|qu|=|\bar q\bar u|$ and $|\uparrow_q^u\xi|\leqslant\epsilon$.
Consequently, we can draw two triangles $\triangle \hat p\hat q\hat u$ and $\triangle \hat s\hat q\hat u$ in $\Bbb S_k^2$
with $|\hat p\hat q|=|pq|$, $|\hat u\hat q|=|uq|$, $|\hat s\hat q|=|sq|$, $\angle\hat p\hat q\hat u=|\uparrow_q^p\uparrow_q^u|$, and $\angle\hat s\hat q\hat u=|\uparrow_q^s\uparrow_q^u|$.
Note that $|\angle\hat p\hat q\hat u-\angle\bar p\bar q\bar u|\leqslant\epsilon$ and $|\angle\hat s\hat q\hat u-\angle\bar s\bar q\bar u|\leqslant\epsilon$. Then via the Law of Cosines on $\Bbb S_k^2$, there is a positive constant $C$ such that
$$\limsup_{s\to q}\left|\frac{|\hat p\hat u|+|\hat u\hat s|-(|\bar p\bar u|+|\bar u\bar s|)}{|qs|}\right|\leqslant C\epsilon,$$
and thus by letting $\epsilon\to 0$ we can conclude that
$$|\hat p\hat u|+|\hat u\hat s|=|\bar p\bar u|+|\bar u\bar s|+o(|qs|)=|\bar p\bar s|+o(|qs|).$$
Moreover, by (4-2) on $\bar X$ and $X$ we have that
$$|pu|\leqslant |\hat p\hat u| \text{ and } |us|\leqslant |\hat u\hat s|\eqno{(5.3)}$$
(note that $p,q, u\in \bar X$ and $s,q,u\in X$).
Then by the definition of $D(X)$,
$$|ps|\leqslant |pu|+|us|\leqslant |\hat p\hat u|+|\hat u\hat s|=|\bar p\bar s|+o(|qs|).$$

We now verify Claim 2. Similarly, we consider only the case where $p\in \bar X^\circ$
(if $p\in\partial X$, (5.2) is true by (4-2) on $X$), and we can assume that $0<\angle pqr<\pi$.
Then $[pq]\cap\partial X$ is a single point $v$ (by Fact 5.2), and for the triangle $\triangle \tilde q\tilde v\tilde s\subset\Bbb S_k^2$ with $|\tilde q\tilde v|=|qv|$, $|\tilde q\tilde s|=|qs|$ and $|\tilde v\tilde s|=|vs|$, we have that
$$0<\angle\tilde q\tilde v\tilde s\leqslant|\uparrow_v^q\uparrow_v^s|\ \text{ and }\ 0<\angle\tilde v\tilde q\tilde s\leqslant|\uparrow_q^v\uparrow_q^s|\ (=\angle pqr<\pi)\eqno{(5.4)}$$
by (4-2) on $X$ (note that $[vq], [qs]$ and any $[vs]$ belong to $X$ by Fact 5.2, and there is a unique minimal geodesic between $q$ and $v$ by the last possible case before Claim 1). We prolong $[\tilde q \tilde v]$ to $\tilde p$ with $\angle\tilde s\tilde v\tilde p=\pi-\angle\tilde q\tilde v\tilde s$ and $|\tilde v\tilde p|=|vp|$. It is clear that $|\tilde p\tilde s|\leqslant |\bar p\bar s|$ by the second part of (5.4). Then, for $s\in [qr]$ sufficiently close to $q$, it suffices to show that
$$|ps|\leqslant |\tilde p\tilde s|+o(|qs|^2).\eqno{(5.5)}$$
\hskip4.5mm According to the last possible case before Claim 1, we first note that $$|\uparrow_v^q\uparrow_v^s|+|\uparrow_v^s\uparrow_v^p|=\pi \text{ (due to $|\uparrow_v^q\uparrow_v^p|=\pi$ in $D(\Sigma_vX))$,}$$
and $\delta\triangleq|\uparrow_v^q\uparrow_v^s|\to 0$ as $s\to q$ (because there is a unique minimal geodesic between $v$ and $q$). And note that $\uparrow_v^s\in (\Sigma_v X)^\circ$ (for $s$ close to $q$) and $\uparrow_v^p\in (\Sigma_v \bar X)^\circ$, so there is $\xi\in\partial\Sigma_vX$ such that
$|\uparrow_v^s\xi|+|\xi\uparrow_v^p|=|\uparrow_v^s\uparrow_v^p|=\pi-\delta$.
Again by Fact 5.3, as $\delta\to 0$, there is
$\epsilon(\delta)\to 0$ and $u\in\partial X$ with $|vu|\to 0$ such that $|\uparrow_v^u\xi|<\epsilon(\delta)$.
Here, a key observation is: since $D(\Sigma_vX)$ is a spherical suspension with diameter equal to $\pi$, from `$|\uparrow_v^s\xi|+|\xi\uparrow_v^p|=\pi-\delta$' and `$|\uparrow_v^u\xi|<\epsilon(\delta)$' we can derive that
$$|\uparrow_v^s\uparrow_v^u|+|\uparrow_v^u\uparrow_v^p|<|\uparrow_v^u\uparrow_v^p|+|\uparrow_v^u\uparrow_v^q|-\delta+ o(\delta)=\pi-\delta+ o(\delta).$$
Consequently, we have a further estimate as follows:
$$|\uparrow_v^s\uparrow_v^u|+|\uparrow_v^u\uparrow_v^p|<\pi-\angle\tilde q\tilde v\tilde s+ o(|qs|)=\angle\tilde s\tilde v\tilde p+ o(|qs|).\eqno{(5.6)}$$
In fact, note that $\angle\tilde q\tilde v\tilde s\leqslant |\uparrow_v^q\uparrow_v^s|=\delta$ (see the first part of (5.4)) with $\delta\to0$ as $|qs|\to 0$, and it is easy to see that $\angle\tilde q\tilde v\tilde s=O(|qs|)$, a positive infinitesimal of the same order as $|qs|$. These enable us to have another key observation that $\angle\tilde q\tilde v\tilde s-\delta+o(\delta)\leqslant o(|qs|)$, which implies (5.6).

We next let $\tilde u$ be the point in $[\tilde p\tilde s]$ such that $\angle \tilde p\tilde v\tilde u=|\uparrow_v^p\xi|$.
Then for $s\in[qr]$ sufficiently close to $q$, it also holds that $|\tilde v\tilde u|=O(|qs|)$, and thus there is $u\in\partial X$
such that $|vu|=|\tilde v\tilde u|$ and $|\uparrow_v^u\xi|<\epsilon(\delta)$. We consider $\triangle \tilde p\tilde v\hat u$ and $\triangle \tilde s\tilde v\hat u$ in $\Bbb S_k^2$ with  $|\hat u\tilde v|=|uv|$, $\angle\tilde p\tilde v\hat u=|\uparrow_v^p\uparrow_v^u|$, and
$\angle\tilde s\tilde v\hat u=\angle\tilde s\tilde v\tilde p-\angle\tilde p\tilde v\hat u$ which implies $|\uparrow_v^s\uparrow_v^u|\leqslant \angle\tilde s\tilde v\tilde p-\angle\tilde p\tilde v\hat u+o(|qs|)=\angle\tilde s\tilde v\hat u+o(|qs|)$ (by (5.6)). Similar to (5.3), we have that
$$|pu|\leqslant |\tilde p\hat u| \text{ and } |us|\leqslant |\hat u\tilde s|+o(|qs|^2).$$
Moreover, note that
$\angle\tilde u\tilde v\hat u=|\angle\tilde p\tilde v\tilde u-\angle\tilde p\tilde v\hat u|=||\uparrow_v^p\xi|-|\uparrow_v^p\uparrow_v^u||\leqslant |\uparrow_v^u\xi|<\epsilon(\delta)$,
and $|\tilde u\hat u|$ as an infinitesimal of $|qs|$ has the same order as $|\tilde v\tilde u|\cdot\sin\angle\tilde u\tilde v\hat u$ (note that $|\tilde v\tilde u|=|\tilde v\hat u|=O(|qs|)$, and $\epsilon(\delta)\to 0$ as $|qs|\to0$). It follows that
$$|\tilde u\hat u|=o(|qs|).$$
This enables us to see that
$|\tilde p\hat u|+|\hat u\tilde s|=|\tilde p\tilde s|+o(|qs|^2)$ in $\triangle \tilde p\tilde s\hat u$ (a very narrow triangle). Then, by the definition of $D(X)$,
$$|ps|\leqslant |pu|+|us|\leqslant |\tilde p\hat u|+|\hat u\tilde s|+o(|qs|^2)=|\tilde p\tilde s|+o(|qs|^2).$$
I.e., (5.5) is verified, and thus the proof is completed.
\end{proof}

%%%%%%%%%%%%%%%%%%%%%%%%%%%%%%%%%%%%%%%%%%%%%%%%%% Section 6%%%%%%%%%%%%%%%%%%%%%%%%%%%%%%%%%%%%%%%%%%%%%%%%%%%%%%%%%%%

\section{Proofs for Globalization Theorems}

This section aims to supply new proofs for the Globalization Theorem only under `complete' and `geodesic' respectively. The new proofs are inspired by the idea of the proofs of Theorems A and B where `locally complete' and `locally geodesic' play crucial role respectively.

\subsection{Only under `complete'}

In this subsection, we will re-prove the Globalization Theorem only under the `complete' condition:

\begin{theorem}[{[BGP]}]\label{main2}
Let $X$ be a complete Alexandrov space with curvature $\geqslant k$. Then for any four distinct points $a,b,c,d\in X$ \footnote{If $k>0$, we  need only consider the case where $|pq|<\frac\pi{\sqrt k}$ for any $p,q\in  X$; and once the proof is completed, we can prove that $|pq|\leqslant\frac\pi{\sqrt k}$ ([BGP]).}, we have that
$\widetilde{\angle}_kbac+\widetilde{\angle}_kbad+\widetilde{\angle}_kcad\leqslant 2\pi.$
\end{theorem}

For convenience, in this subsection, $X$ always denotes the space in Theorem 6.1. Each known proof of Theorem 6.1 considers first the ideal case where $X$ is in addition geodesic. Recently, a new proof of Theorem 6.1 for the ideal case appears in [HSW]. Inspired by this new proof and the proof of Theorem A, here we present a direct proof for Theorem 6.1 without `geodesic'.
For the purpose, we need further to develop `imaginary' angles, and introduce `imaginary' distance functions.

\subsubsection{`Imaginary' angles and distance functions}

\begin{defn}\label{angle}{\rm  For three distinct points $p,q,r\in X$, an {\it imaginary angle} at $p$ is defined to be
$\measuredangle qpr\triangleq\lim\limits_{x,y\to 0}\omega_k[p_q^r](x,y)$ (if this limit exists), where $$\omega_k[p_q^r](x,y)\triangleq\inf\{\lim_{i\to\infty}\widetilde{\angle}_ku_ipv_i|\text{ $\mathcal{E}_x^{pq}(u_i)\to0$ and $\mathcal{E}_y^{pr}(v_i)\to 0$ as $i\to\infty$}\}\leqslant\pi$$
with $x\in(0,|pq|]$ and $y\in(0,|pr|]$ (refer to Section 2 for the `error' function $\mathcal{E}$).}
\end{defn}

\begin{defn}\label{chord}{\rm For three distinct points $p,q,r\in X$, an {\it imaginary distance function from $p$  to $(q,r)$} with respect to $d\in[0,|qr|]$  is defined to be
$$	\dist[_p^{qr}](d)\triangleq\inf\{\lim_{i\to\infty}|ps_i||\text{ $\mathcal{E}_d^{qr}(s_i)\to 0$ as $i\to\infty$}\}.$$
Moreover, we let $\widetilde{\dist}_k[_p^{qr}](\cdot)\triangleq\dist[_{\tilde{p}}^{\tilde{q}\tilde{r}}](\cdot)$, where $\tilde{p},\tilde{q},\tilde{r}\in\mathbb{S}_k^2$
with $|\tilde p\tilde q|=|pq|, |\tilde p\tilde r|=|pr|, |\tilde q\tilde r|=|qr|$. It is not hard to check that $\dist[_p^{qr}](\cdot)$ is 1-Lipschitz on $[0,|qr|]$. }
\end{defn}

These concepts make it possible for us to argue as on geodesic spaces. And, similar to $R_p(q)$ in Section 2, we will use $R(p)$ to denote the radius of good neighborhood around $p$, i.e.
$$R(p)\triangleq\max\{R|\text{ $\widetilde{\angle}_kbac+\widetilde{\angle}_kbad+\widetilde{\angle}_kcad\leqslant 2\pi$ for any distinct points $a,b,c,d\in B(p,R)$}\}.$$

Similar to some basic properties of angles in the case where $X$ is in addition geodesic (cf. [BGP]),
we have the following four propositions on the `imaginary' angles and distance functions.

\begin{prop}\label{equiv1} Given three distinct points $p,q,r\in X$, we have that $\omega_k[p_q^r](x_1,y)\geqslant\omega_k[p_q^r](x_2,y)$ for any $0<x_1<x_2<R(p)$ and $0<y<R(p)$.
As a result,

{\rm (6.4.1)} $\measuredangle qpr$ in Definition 6.2 is well-defined, and $\measuredangle qpr\geqslant\widetilde{\angle}_kqpr$ if $q,r\in B(p,R(p))$;

{\rm (6.4.2)} $\measuredangle bac+\measuredangle bad+\measuredangle cad\leqslant 2\pi$ for all four distinct points $a,b,c,d\in X$.
\end{prop}

In the special case where $X$ is in addition of finite dimension, $\measuredangle qpr$ is just the minimum of angles between each minimal geodesic $[pq]$ and each $[pr]$.

\begin{proof}[Proof] By Definition 6.2, there is $\{u_{1i}, v_i\}_{i=1}^\infty\subset X$ such that  $\mathcal{E}_{x_1}^{pq}(u_{1i}),\mathcal{E}_y^{pr}(v_i)\to 0$ and $\widetilde{\angle}_ku_{1i}pv_i\to\omega_k[p_q^r](x_1,y)$ as $i\to\infty$. Then we select $\{u_{2i}\}_{i=1}^\infty$ such that $\mathcal{E}_{x_2-x_1}^{u_{1i}q}(u_{2i})\to 0$ and thus $\mathcal{E}_{x_2}^{pq}(u_{2i})\to 0$ as $i\to\infty$  (by (2.1)). Applying (0.1) to $\{u_{1i}; p,u_{2i},v_i\}$, we have that $\limsup\limits_{i\to\infty}\left(\widetilde{\angle}_kv_iu_{1i}p+\widetilde{\angle}_kv_iu_{1i}u_{2i}\right)\leqslant\pi$; then by Lemma 2.1 we can conclude that
\begin{align*}
	\omega_k[p_q^r](x_1,y)=\lim_{i\to\infty}\widetilde{\angle}_kv_ipu_{1i}\geqslant\limsup_{i\to\infty}\widetilde{\angle}_kv_ipu_{2i}\geqslant\omega_k[p_q^r](x_2,y).
\end{align*}

Similarly, $\omega_k[p_q^r](x,y)$ is also decreasing with respect to $y$ when $0<x, y<R(p)$. This implies that $\measuredangle qpr$ is well-defined,
and $\measuredangle qpr\geqslant\widetilde{\angle}_kqpr$ if $q,r\in B(p,R(p))$. And by taking into account (0.1), we can conclude that $\measuredangle bac+\measuredangle bad+\measuredangle cad\leqslant 2\pi$.
\end{proof}

\begin{prop}\label{1var} Given three distinct points $p,q,r\in X$, we have that
	\[
	\dist[_p^{qr}](t)\leqslant|pq|-\cos\measuredangle pqr\cdot t+o(t) \text{ with } t\in(0,|qr|),
	\]
	where $o(t)/t\to0$ as $t\to 0^+$ and $o(t)$ depends only on $R(q)$.
\end{prop}

\begin{proof}[Proof]
	Let $l\triangleq\min\{|pq|,\frac{R(q)}{2}\}$, and consider only $t\ll l$. By Definition 6.2, there is $\{p_i, s_i\}_{i=1}^\infty\subset X$ such that $\mathcal{E}_l^{qp}(p_i), \mathcal{E}_t^{qr}(s_i)\to 0$ and $\widetilde{\angle}_kp_iqs_i\to\omega_k[q_p^r](l,t)$ as $i\to\infty$. Then by Definition 6.3, (6.4.1) and the first variation formula on $\mathbb{S}_k^2$, we have that
	\begin{align*}
		\dist[_p^{qr}](t)&\leqslant\liminf_{i\to\infty}|ps_i|\leqslant\liminf_{i\to\infty}\left(|pp_i|+|p_is_i|\right)\\
		&\leqslant\lim_{i\to\infty}\left(|pp_i|+|p_iq|-\cos\measuredangle pqr\cdot|qs_i|+o(|qs_i|)\right)\\
		&=|pq|-\cos\measuredangle pqr\cdot t+o(t)
	\end{align*}
(note that $o(t)$ can be chosen to depend only on  $R(q)$).
\end{proof}

\begin{prop}\label{minimum}
	Given three distinct points $p,q,r\in X$, and $d_1\in(0,|pq|]$ and $d_2\in(0,|pr|]$, we have that
$\measuredangle qpr\leqslant\inf\left\{\lim\limits_{i\to\infty}\measuredangle u_ipv_i\left|\text{ $\mathcal{E}_{d_1}^{pq}(u_i),\mathcal{E}_{d_2}^{pr}(v_i)\to 0$ as $i\to\infty$}\right. \right\}$.
\end{prop}

\begin{proof}[Proof] Assume that $\{u_i,v_i\}_{i=1}^\infty\subset X$ satisfies that $\mathcal{E}_{d_1}^{pq}(u_i),\mathcal{E}_{d_2}^{pr}(v_i)\to 0$ as $i\to\infty$ and $\lim\limits_{i\to\infty}\measuredangle u_ipv_i$ exists.
Then, by Proposition \ref{equiv1}, it suffices to show that $\omega_k[p_q^r](\delta,\delta)\leqslant\lim\limits_{i\to\infty}\measuredangle u_ipv_i$ for any sufficiently small $\delta>0$.
Note that by Definition 6.2, there is $\{s_i, t_i\}_{i=1}^\infty\subset X$ such that $\mathcal{E}_{\delta}^{pu_i}(s_i),\mathcal{E}_{\delta}^{pv_i}(t_i)\to 0$ (and thus $\mathcal{E}_{\delta}^{pq}(s_i),\mathcal{E}_{\delta}^{pr}(t_i)\to 0$) and $\widetilde{\angle}_ks_ipt_i-\omega_k[p_{u_i}^{v_i}](\delta,\delta)\to 0$  as $i\to\infty$. By Proposition \ref{equiv1}, we have

\vskip2mm

\ \hskip2cm	$\omega_k[p_q^r](\delta,\delta)\leqslant\liminf\limits_{i\to\infty}\widetilde{\angle}_ks_ipt_i=\liminf\limits_{i\to\infty}\omega_k[p_{u_i}^{v_i}](\delta,\delta)\leqslant\lim\limits_{i\to\infty}\measuredangle u_ipv_i.$
\end{proof}

\begin{prop}\label{complementary} For three distinct points $p,q,r\in X$ and $d\in(0,|qr|)$ with $\dist[_p^{qr}](d)>0$, suppose that $\mathcal{E}_d^{qr}(s_i)\to 0$ as $i\to\infty$, and $\{R(s_i)\}_{i=1}^\infty$ has a positive lower bound. Then the followings hold:
	
	{\rm (6.7.1)} $\limsup\limits_{i\to\infty}(\measuredangle ps_iq+\measuredangle ps_ir)\leqslant\pi$.
	
	{\rm (6.7.2)} If in addition $|ps_i|\to\dist[_p^{qr}](d)$ as $i\to\infty$, then
	\[
	-\dist[_p^{qr}]'_{-,\min}(d)\leqslant-\cos\liminf_{i\to\infty}\measuredangle ps_iq\ \text{ and }\ \dist[_p^{qr}]'_{+,\max}(d)\leqslant-\cos\liminf_{i\to\infty}\measuredangle ps_ir;
	\]
	
   {\rm (6.7.3)} as a result,  $\dist[_p^{qr}]'_{+,\max}(d)-\dist[_p^{qr}]'_{-,\min}(d)\leqslant 0$, and ``$=$'' implies
     $$\hskip1cm	-\dist[_p^{qr}]'_{-,\min}(d)=-\cos\lim_{i\to\infty}\measuredangle ps_iq\ \text{ and }\ \dist[_p^{qr}]'_{+,\max}(d)=-\cos\lim_{i\to\infty}\measuredangle ps_ir$$
	\hskip1.5cm(and $\lim\limits_{i\to\infty}\measuredangle  ps_iq+\lim\limits_{i\to\infty}\measuredangle ps_ir=\pi$).
\end{prop}

Note that (6.7.3) is under the assumption of (6.7.2). And note that in the case where $X$ is geodesic, (6.7.1) can be represented by `$\angle psq+\angle psr=\pi$' (see Subsection 1.3).

\begin{proof}[Proof] (6.7.1)\ \ From the conditions, we can assume that $\min\{R(s_i),|ps_i|,|qs_i|,|rs_i|\}_{i=1}^\infty>2c>0$.
By Proposition \ref{minimum}, there is  $\{p_i,q_i,r_i\}_{i=1}^\infty$ such that, as $i\to\infty$, $\mathcal{E}_c^{s_ip}(p_i),\mathcal{E}_c^{s_iq}(q_i),\mathcal{E}_c^{s_ir}(r_i)\to 0$ and $\limsup\limits_{i\to\infty}(\measuredangle ps_iq+\measuredangle ps_ir)\leqslant\limsup\limits_{i\to\infty}(\measuredangle p_is_iq_i+\measuredangle p_is_ir_i)$.
Then, by (6.4.2), it suffices to show $\lim\limits_{i\to\infty}\measuredangle q_is_ir_i=\pi$. Since $q_i,r_i\in B(s_i, R(s_i))$, it is indeed true by (6.4.1) (note that
$\widetilde{\angle}_kq_is_ir_i\to\pi$ by (2.1)).

\vskip1mm
	
(6.7.2)\ \ For any small $\tau>0$, by Definition 6.3 and (2.1) we have that $\dist[_p^{qr}](d+\tau)\leqslant\liminf\limits_{i\to\infty}\dist[_p^{s_ir}](\tau)$.
And by Proposition \ref{1var}, there is an $o(\tau)$ depending only on $\min\{R(s_i)\}_{i=1}^\infty>0$ such that
$$\dist[_p^{s_ir}](\tau)\leqslant |ps_i|-\cos\measuredangle ps_ir\cdot \tau+o(\tau).$$
Since $\lim\limits_{i\to\infty}|ps_i|=\dist[_p^{qr}](d)$, it follows that
$$\dist[_p^{qr}](d+\tau)\leqslant\dist[_p^{qr}](d)-\cos\liminf\limits_{i\to\infty}\measuredangle ps_ir\cdot \tau+o(\tau).$$
I.e., $\dist[_p^{qr}]'_{+,\max}(d)\leqslant-\cos\liminf\limits_{i\to\infty}\measuredangle ps_ir$;
and similarly, $-\dist[_p^{qr}]'_{-,\min}(d)\leqslant-\cos\liminf\limits_{i\to\infty}\measuredangle ps_iq$.

\vskip1mm
	
(6.7.3)\ \ Note that, by (6.7.1),
$$\liminf_{i\to\infty}\measuredangle ps_iq+\liminf_{i\to\infty}\measuredangle ps_ir\leqslant \liminf_{i\to\infty}(\measuredangle ps_iq+\measuredangle ps_ir)\leqslant \limsup_{i\to\infty}(\measuredangle ps_iq+\measuredangle ps_ir)\leqslant\pi.$$  So, it follows from (6.7.2) that $\dist[_p^{qr}]'_{+,\max}(d)-\dist[_p^{qr}]'_{-,\min}(d)\leqslant 0$. And if the equality holds, then both equalities in (6.7.2) hold, and $$\liminf\limits_{i\to\infty}\measuredangle ps_iq+\liminf\limits_{i\to\infty}\measuredangle ps_ir=\pi=\lim\limits_{i\to\infty}(\measuredangle ps_iq+\measuredangle ps_ir).$$
Here the second `$=$' implies that $\liminf\limits_{i\to\infty}\measuredangle ps_iq=\pi-\limsup\limits_{i\to\infty}\measuredangle ps_ir$, which plus the first `$=$' implies that $\lim\limits_{i\to\infty}\measuredangle ps_ir$ exists, so does $\lim\limits_{i\to\infty}\measuredangle ps_iq$. It then follows that $-\dist[_p^{qr}]'_{-,\min}(d)=-\cos\lim\limits_{i\to\infty}\measuredangle ps_iq$  and
$\dist[_p^{qr}]'_{+,\max}(d)=-\cos\lim\limits_{i\to\infty}\measuredangle ps_ir $, and $\lim\limits_{i\to\infty}\measuredangle  ps_iq+\lim\limits_{i\to\infty}\measuredangle ps_ir=\pi$.
\end{proof}

%\begin{remark}\label{comple1}{\rm  In Proposition \ref{complementary}, if the $p,q,r$ are replaced with $p_i,q_i,r_i$ satisfying
%$\liminf\limits_{i\to\infty}|q_ir_i|>d$ and $\liminf\limits_{i\to\infty}\dist[_{p_i}^{q_ir_i}](d)>0$, then from the proof of (6.7.1) we can conclude that it still holds.}
%\end{remark}

We are now ready to prove Theorem \ref{main2} (note that Propositions 6.4-6.7 have nothing to do with `complete').

\subsubsection{Proof of Theorem 6.1}

By (6.4.2), it suffices to show that any `imaginary' $\measuredangle pqr$ satisfies $\measuredangle pqr\geqslant\widetilde{\angle}_kpqr$. We will argue by contradiction, and say that $\measuredangle pqr$ is {\it bad} if
$\measuredangle pqr<\widetilde{\angle}_kpqr$; and denote by
\[
	S_p\triangleq\{o\in X|\ \forall\ \delta>0,\exists\ r_1,r_2\in B(o, \delta)\mbox{ s.t. }\measuredangle pr_1r_2\mbox{ is bad}\}.
\]
Then inspired by the idea to look for the $\bar s$ in the proof of Lemma 2.2 and the idea of Corollary 2.2 in [HSW], we have the following observation.
For the sake of simplicity and due to the similarity, we consider only the case where $k=0$ in Theorem 6.1, and let $\widetilde \angle pqr$ denote $\widetilde\angle_k pqr$.

\begin{lemma}\label{core2.0}
Given three distinct points $p,r_1,r_2\in X$ with $|r_1r_2|<\max\{R(r_1),R(r_2)\}$,  if $\measuredangle pr_1r_2$ is bad, then there is $s\in X$ with $|r_1s|+|r_2s|$ arbitrarily close to $|r_1r_2|$ such that $s\in S_p$ and $|ps|\leqslant\max\{|pr_1|,|pr_2|\}$.
\end{lemma}

\begin{proof}[Proof] Note that the badness of $\measuredangle pr_1r_2$ implies $\rho_0\triangleq|pr_1|+|pr_2|-|r_1r_2|>0$.
Then for any $0<\delta\ll \rho_0$, by Lemma \ref{core2} below there is $d_1\in(0, |r_1r_2|)$ and $s_1\in X$ with $\mathcal{E}_{d_{1}}^{r_1r_2}(s_{1})<\frac12\delta$
and $d_1<|pr_1|$ or $|r_1r_2|-d_1<|pr_2|$, say $|r_1r_2|-d_1<|pr_2|$, such that $\measuredangle ps_{1}r_2$ is bad. Moreover, $|ps_1|<\max\{|pr_1|,|pr_2|\}$, and $s_1$ can be selected such that $|s_1r_2|<|pr_2|$ and $|s_1r_2|<\max\{R(s_1),R(r_2)\}$.
Then by Lemma \ref{core2} again, there is $d_2\in(0, |s_1r_2|)$ and $s_2$ with $\mathcal{E}_{d_{2}}^{s_1r_2}(s_{2})<\frac14\delta$,
$|ps_2|<\max\{|pr_1|,|pr_2|\}$, and $|s_2r_2|<|pr_2|$ and $|s_2r_2|<\max\{R(s_2),R(r_2)\}$ such that $\measuredangle ps_{2}r_2$ is bad.
If $d_1+d_2\geqslant|r_1r_2|$ for all such $s_1$ and $s_2$, then $s_2$ converges to $r_2$ as $\delta\to0$, i.e. $r_2\in S_p$; if $d_1+d_2<|r_1r_2|$, then we can similarly locate an $s_3$. Step by step, if $r_2\not\in S_p$, we can locate $\{s_i\}_{i=1}^\infty$ and $\{d_i\}_{i=1}^\infty$ such that $\mathcal{E}_{d_{i+1}}^{s_ir_2}(s_{i+1})<\frac1{2^{i+1}}\delta$, $\sum\limits_{i=1}^\infty d_i\leqslant|r_1r_2|$, $|ps_i|<\max\{|pr_1|,|pr_2|\}$, and $\measuredangle ps_{i}r_2$ is bad.
Note that  $\sum\limits_{i=1}^\infty |s_is_{i+1}|<\sum\limits_{i=1}^\infty \left(d_i+\frac\delta{2^i}\right)\leqslant|r_1r_2|+\delta$,
which implies that $\{s_i\}_{i=1}^\infty$ is a Cauchy sequence and thus $s_i$ converges to an $s$ as $i\to\infty$ (because $X$ is complete).
And note that $\mathcal{E}_{d_1+\cdots+d_{i}}^{r_1r_2}(s_{i})<\delta$ by (2.1), and so $|s_is_{i+1}|<\frac12\rho_0-2\delta<|ps_i|$ for large $i$. This implies that $\measuredangle ps_{i+1}s_{i}$ is also bad by Lemma \ref{core2}, i.e. $s$ is just a desired point.
\end{proof}

\begin{lemma}\label{core2}
Given three distinct points $p,r_1,r_2\in X$ with $|r_1r_2|<\max\{R(r_1),R(r_2)\}$, if $\measuredangle pr_1r_2$ is bad, then there is $d\in(0,|r_1r_2|)$ and $\{s_i\}_{i=1}^\infty$ with $\mathcal{E}_d^{r_1r_2}(s_i)\to0$ as $i\to\infty$ such that each $\measuredangle ps_ir_1$ is bad or each $\measuredangle ps_ir_2$ is bad, and $|ps_i|<\max\{|pr_1|,|pr_2|\}$ \footnote{For $k>0$, it needs to be modified to $|ps|\leqslant\max\{|pr_1|,|pr_2|\}+o(|r_1r_2|)$, where $o(|r_1r_2|)$ depends only on $\max\{|pr_1|,|pr_2|\}<\pi$ (see footnote 7). \label{foot2}}. In particular, if $d\leqslant|r_1p|$ \footnote{For $k>0$, it needs to be modified to $d\leqslant\min\{|r_1p|,\frac{\pi}{\sqrt{k}}-|r_1p|\}$.\label{foot1}}, then each $\measuredangle ps_ir_1$  is bad (and if $|r_1r_2|-d\leqslant|r_2p|$, then each $\measuredangle ps_ir_2$ is bad).
\end{lemma}

The lemma is an imitation of Lemma 6.12 below.

\begin{proof}[Proof]
By Proposition \ref{1var} and the first variation formula on $\mathbb{R}^2$, the badness of $\measuredangle pr_1r_2$ implies that the function $\dist[_p^{r_1r_2}]-\widetilde{\dist}[_p^{r_1r_2}]$ attains a negative local minimum at some $d\in(0,|r_1r_2|)$. Let $\triangle\tilde{p}\tilde{r}_1\tilde{r}_2\subset \Bbb R^2$ satisfy
$|\tilde p\tilde r_1|=|pr_1|, |\tilde p\tilde r_2|=|pr_2|$ and $|\tilde r_1\tilde r_2|=|r_1r_2|$, and $\tilde{s}\in[\tilde{r}_1\tilde{r}_2]$ satisfy $|\tilde{r}_1\tilde{s}|=d$. Then the minimality at $d$ implies that
\begin{align}
\dist[_p^{r_1r_2}]'_{+,\min}(d)\geqslant-\cos\angle\tilde{p}\tilde{s}\tilde{r}_2\ \text{ and }  -\dist[_p^{r_1r_2}]'_{-,\max}(d)\geqslant-\cos\angle\tilde{p}\tilde{s}\tilde{r}_1.\label{1}
\end{align}
By Definition 6.3, there is $\{s_i\}_{i=1}^\infty$ such that $\mathcal{E}_d^{r_1r_2}(s_i)\to 0$ and $|ps_i|\to\dist[_p^{r_1r_2}](d)$ as $i\to\infty$ (but any such $\{s_i\}_{i=1}^\infty$ might NOT contain a Cauchy sequence). Since $|r_1r_2|<\max\{R(r_1),R(r_2)\}$, we can assume that $\{R(s_i)\}_{i=1}^\infty$ have a positive lower bound.
Then by the first part of (6.7.3), (6.1) is actually
\begin{align*}
\dist[_p^{r_1r_2}]'_{+}(d)=-\cos\angle\tilde{p}\tilde{s}\tilde{r}_2\ \text{ and }  -\dist[_p^{r_1r_2}]'_{-}(d)=-\cos\angle\tilde{p}\tilde{s}\tilde{r}_1
\end{align*}
(note that $\cos\angle\tilde{p}\tilde{s}\tilde{r}_1+\cos\angle\tilde{p}\tilde{s}\tilde{r}_2=0$);
so, by the case where ``$=$'' holds in (6.7.3),
 $$	\lim_{i\to\infty}\measuredangle ps_ir_2=\angle\tilde{p}\tilde{s}\tilde{r}_2\ \text{ and }\
 \lim_{i\to\infty}\measuredangle ps_ir_1=\angle\tilde{p}\tilde{s}\tilde{r}_1. \eqno{(6.2)}$$

On the other hand, since $\dist[_p^{r_1r_2}](d)<|\tilde p\tilde s|$,
there is $\tilde p_i\in [\tilde p\tilde s]$ for large $i$ such that $|\tilde p_i\tilde s|=|ps_i|\ (<|\tilde p\tilde s|)$ and $\tilde p_i\to \tilde p'$ with $|\tilde p'\tilde s|=\dist[_p^{r_1r_2}](d)$; and so we can assume that $|\tilde p_i\tilde r_1|<| pr_1|$ (note that $|\tilde p'\tilde r_1|+|\tilde p'\tilde r_2|<|pr_1|+|pr_2|$), especially when $d\leqslant |pr_1|$. It follows that $|ps_i|<|\tilde p\tilde s|<\max\{|pr_1|,|pr_2|\}$, and $\widetilde{\angle}ps_ir_1>\angle\tilde{p}\tilde{s}\tilde{r}_1$ which plus (6.2) implies that $\measuredangle ps_i{r}_1$ is bad. This completes the proof.
\end{proof}

\begin{remark}\label{6.11}{\rm  In the proof right above, the condition `$|r_1r_2|<\max\{R(r_1),R(r_2)\}$' just guarantees the property that $\{R(s_i)\}_{i=1}^\infty$ have a positive lower bound. So, if the condition is omitted but the property is satisfied, the proof will also go through.}
\end{remark}

Next, we will finish the proof of Theorem 6.1 using Lemmas 6.8 and 6.9, and the rough idea of the proof is almost the same as that of the proof in [HSW].

\renewcommand{\proofname}{Proof of Theorem \ref{main2}}\begin{proof}
By (6.4.2), it suffices to show that, for any three distinct points $p,q,r\in X$, $\measuredangle pqr\geqslant\widetilde{\angle}pqr$.
If some $\measuredangle pqr$ is bad, we will get a contradiction through the following two steps.

\vspace{0.1cm}
\noindent {\bf Step I}: To locate an $s\in S_p$ due to the badness of $\measuredangle pqr$.
\vspace{0.1cm}

Lemma 6.8 implies the existence of a desired $s$  if $|qr|<R(q)$ or $R(r)$ \footnote{If $X$ is in addition geodesic, $s$ can be located in $[qr]$ directly without the need of `$|qr|<R(q)$ or $R(r)$' ([HSW]).}. In other cases, let $R_x\triangleq\min\{\frac{1}{2}R(x),|px|\}$ for $x\in X$, and consider the following two cases.
Note that $R_q>0$.

\vspace{0.1cm}
\noindent{\bf(a)} $\dist[_p^{qr}](R_q)\geqslant\widetilde{\dist}[_p^{qr}](R_q)$:
\vspace{0.1cm}

By Proposition \ref{1var}, the badness of $\measuredangle qpr$ implies that the function $\dist[_p^{qr}]-\widetilde{\dist}[_p^{qr}]$ attains a negative local minimum at some $d\in(0,R_q)$.
Note that for any $\{u_i\}_{i=1}^\infty$ with $\lim\limits_{i\to\infty}\mathcal{E}_d^{qr}(u_i)=0$ and $\lim\limits_{i\to\infty}|pu_i|=\dist[_p^{qr}](d)$,  $R(u_i)>\frac12 R(q)>|u_iq|<|pq|$ for large $i$.
Then the proof of Lemma 6.9 also implies that $\measuredangle pu_iq$ is bad (cf. Remark 6.10), and thus we can locate an $s\in S_p$ by Lemma \ref{core2.0}.

\vspace{0.1cm}
\noindent{\bf (b)} $\dist[_p^{qr}](R_q)<\widetilde{\dist}[_p^{qr}](R_q)$:
\vspace{0.1cm}

By Definition 6.3, there is $\{u_i\}_{i=1}^\infty$ such that $\mathcal{E}_{R_q}^{qr}(u_i)\to 0$ and $|pu_i|\to\dist[_p^{qr}](R_q)$ as $i\to\infty$.
Since $\dist[_p^{qr}](R_q)<\widetilde{\dist}[_p^{qr}](R_q)$, by Lemma 2.1 we can conclude that
$$\lim_{i\to\infty}(\widetilde{\angle}pu_iq+\widetilde{\angle}pu_ir)>\pi.$$
On the other hand, since $|qu_i|\to R_q$, $\liminf\limits_{i\to\infty} R(u_i)\geqslant R_q$  and so by (6.7.1) we have
$$\limsup_{i\to\infty}(\measuredangle pu_iq+\measuredangle pu_ir)\leqslant\pi.$$
It follows that $\measuredangle pu_iq$ or $\measuredangle pu_ir$ is bad for large $i$ \footnote{Here we can not get that
$\measuredangle pu_iq$ is bad even if in addition $|qu_i|\ll|qp|$ (cf. Lemma 6.9), so we can not locate the desired $s$ as in Lemma 6.8.}, and one of the following two subcases holds.

{\bf (b1)} $\measuredangle pu_iq$ is bad, or $\measuredangle pu_ir$ is bad and $|qr|-R_q<R(r)$: In this case, we can locate an $s\in S_p$ by Lemma \ref{core2.0} (note that $|u_iq|<R(q)$, and `$|qr|-R_q<R(r)$' implies $|u_ir|<R(r)$).

{\bf (b2)} $\measuredangle pu_ir$ is bad and $|qr|-R_q\geqslant R(r)$: In this case, take a positive $\delta\ll \rho_0\triangleq |pq|+|pr|-|qr|$
(note that $\rho_0>0$ by the badness of $\measuredangle pqr$, otherwise, $\measuredangle pqr=\widetilde\angle pqr=0$).
For large $i$, we can assume that $\mathcal{E}_{R_q}^{qr}(u_i)<\frac12\delta$ and $|u_ir|\geqslant R(r)$,
and we set $s_1\triangleq q$ and $s_2\triangleq u_i$.

\vskip1mm

Note that Step I is finished except when subcase (b2) happens. In the exceptional case, we can repeat the above for $p, s_2, r$.
If we can not find any $s\in S_p$ in this way, then step by step and as in the proof of Lemma 2.2 we can locate $\{s_j\}_{j=1}^\infty$ with $|s_jr|\geqslant R(r)$, $\mathcal{E}_{R_{s_j}}^{s_jr}(s_{j+1})<\frac{\delta}{2^j}$, and $\sum\limits_{j=1}^\infty R_{s_j}<|qr|$ which implies $R_{s_j}\to0$  as $j\to\infty$ and so
$R(s_j)\to0$ because $|ps_{j}|>\frac12\rho_0-2\delta>0$ (note that `$\mathcal{E}_{R_{s_j}}^{s_jr}(s_{j+1})<\frac{\delta}{2^j}$ and $\sum\limits_{j=1}^\infty R_{s_j}<|qr|$' implies that $\mathcal{E}_{R_{s_1}+\cdots+R_{s_j}}^{qr}(s_{j+1})<\delta$ by (2.1)).
However, similar to (2.5), we can say that $\{s_j\}$ is a Cauchy sequence, so $s_j$ has to converge to a point $\bar{s}$ because $X$ is complete. It follows that
$R(s_j)\to R(\bar{s})>0$ as $j\to\infty$; a contradiction.

\vspace{0.1cm}
\noindent {\bf Step II}: To show that the non-emptiness of $S_p$ contradicts `curvature $\geqslant k$' around $p$.

For all $x\in X$, let $\delta(x)\triangleq\min\{\frac{|px|}{2},R(x)\}$. As in [HSW], we claim that
$$\forall\ o\in S_p,\ \exists\ o^\prime\in S_p\cap B(o,\delta(o))\text{ s.t. } |po^\prime|\leqslant|po|-\frac{\delta(o)}{4}. $$
By assuming the claim, we can inductively either locate an $\bar o\in S_p$ with $|p\bar o|<R(p)$ which contradicts `curvature $\geqslant k$' around $p$ (see (6.4.1)), or obtain $\{o_i\}_{i=1}^\infty\subset S_p$ such that
$$R(p)\leqslant |po_{i+1}|\leqslant|po_i|-\frac{\delta(o_i)}{4}\text{ and }|o_io_{i+1}|<\delta(o_i).$$
Note that $|po_{i+1}|\leqslant|po_1|-\sum\limits_{j=1}^i\frac{\delta(o_j)}{4}$ and thus $\sum\limits_{j=1}^\infty\frac{\delta(o_j)}{4}\leqslant |po_1|$, which implies $\lim\limits_{i\to\infty}\delta(o_i)=0$ and $\{o_i\}_{i=1}^\infty$ is a Cauchy sequence. Since $X$ is complete, $o_i$ converges to an $\bar o$ as $i\to\infty$, which implies $\delta(o_i)\to\delta(\bar o)>0$, a contradiction.

It remains to verify the claim above. Set $d=|po|-\frac{\delta(o)}{4}-\epsilon$ with $0<\epsilon\ll \delta(o)$, and fix $r_1$ and $r_2$ with $|or_i|\ll\delta(o)$ and $\measuredangle pr_1r_2$ bad. Then we can locate a desired $o'$ as follows.

\noindent$\bullet$ There is $\{\bar{r}_{1i}\}_{i=1}^\infty$ such that $\lim\limits_{i\to\infty}\mathcal{E}_d^{pr_1}(\bar{r}_{1i})=0$, $\lim\limits_{i\to\infty}|r_2\bar{r}_{1i}|$ exists, and  $\lim\limits_{i\to\infty}\widetilde{\angle}\bar{r}_{1i}r_1r_2=\omega[{r_1}_p^{r_2}](|pr_1|-d,|r_1r_2|)$ (see Definition 6.2). By the badness of $\measuredangle pr_1r_2$ and Proposition \ref{equiv1}, $\lim\limits_{i\to\infty}\widetilde{\angle}\bar{r}_{1i}r_1r_2<\widetilde{\angle}pr_1r_2$.
Note that $\{R(\bar r_{1i})\}_{i=1}^\infty$ has a positive lower bound. Then by Lemma 2.1 plus (6.7.1) and `$\measuredangle r_1\bar{r}_{1i}r_2\geqslant\widetilde{\angle}_kr_1\bar{r}_{1i}r_2$ by (6.4.1)', we derive that $\measuredangle p\bar{r}_{1i}r_2$ is bad for large $i$. Denote by $\bar{r}_1$ such an $\bar{r}_{1i}$. Note that $|po|-\frac{\delta(o)}{4}\geqslant|p\bar{r}_{1}|\ (\stackrel{i\to\infty}{\longrightarrow}d)$, and $|\bar{r}_1o|<\frac{\delta(o)}{4}+2\epsilon$ which implies $|\bar{r}_1r_2|<R(r_2)$.

\noindent$\bullet$ Via Lemma \ref{core2} on $\{p, \bar{r}_1, r_2\}$, there is $\rho\in(0,|\bar{r}_1r_2|)$ and $\hat{s}$ such that $\mathcal{E}_\rho^{\bar{r}_1r_2}(\hat{s})\ll |\bar{r}_1r_2|<\delta(o)$ and $\measuredangle p\hat{s}\bar{r}_1$ is bad (note that $|\bar{r}_1r_2|<|p\bar r_1|$). Note that $\mathcal{E}_{d+\rho}^{po}(\hat{s})\ll \delta(o)$ because $|or_i|\ll \delta(o)$, $\mathcal{E}_d^{pr_1}(\bar{r}_1)\ll \delta(o)$ and $\mathcal{E}_\rho^{\bar{r}_1r_2}(\hat{s})\ll \delta(o)$; and thus we can assume that $|\hat so|<\frac{\delta(o)}{4}+2\epsilon$.

\noindent$\bullet$ If $|p\hat{s}|\leqslant |po|-\frac{\delta(o)}{4}$, then let $\bar{r}_2=\hat{s}$; otherwise, similar to $\bar{r}_1$, we can locate a $\bar{r}_2$ from $\{p, \hat s, \bar r_1\}$ such that $\measuredangle p\bar{r}_2\bar{r}_1$ is bad, $|p\bar{r}_{2}|\leqslant|po|-\frac{\delta(o)}{4}$, $|\bar{r}_2o|<\frac{\delta(o)}{4}+2\epsilon$ (note that $\mathcal{E}_d^{po}(\bar{r}_{1}), \mathcal{E}_{d+\rho}^{po}(\hat{s})\ll \delta(o)$), and $|\bar{r}_2\bar r_1|<R(\bar r_2)$ (note that $|\bar{r}_2\bar r_1|$ is at most about $\frac{\delta(o)}{2}$, while $R(\bar r_2)$ is at least about $\frac{3\delta(o)}{4}$).

\noindent$\bullet$ By applying Lemma \ref{core2.0} to $\{p, \bar{r}_2, \bar{r}_1\}$, we can obtain an $o^{\prime}\in S_p$ satisfying the claim above.
\end{proof}

\renewcommand{\proofname}{Proof}

\subsection{Only under `geodesic'}

In this subsection, we will re-prove the Globalization Theorem only under the `geodesic' condition, i.e. Theorem 6.11.
From its original proof ([Petr2]), our proof is quite different because our proof mainly relies on the new ingredients, Lemma 3.4 and the idea of the proof in [HSW].

\begin{theorem}[{[Petr2]}]\label{6.11}
Let $X$ be a geodesic Alexandrov space with curvature $\geqslant k$. Then for any four distinct points $a,b,c,d\in X$ \footnote{If $k>0$, we also need only consider the case where $|pq|<\frac\pi{\sqrt k}$ for any $p,q\in  X$.}, we have that
$\widetilde{\angle}_kbac+\widetilde{\angle}_kbad+\widetilde{\angle}_kcad\leqslant 2\pi.$
\end{theorem}

For convenience, in this subsection, $X$ always denotes the space in Theorem 6.11. Note that, between any two geodesics $[pq]$ and $[qr]$, we can define a `visible' angle, denoted by $\angle pqr$ (see Subsection 1.3). And for any $[ab], [ac], [ad]\subset X$, it is easy to see that $\angle bac+\angle bad+\angle cad\leqslant 2\pi$
(cf. Proposition 6.4).
Then similarly, to prove Theorem 6.11, we just need to show $\angle pqr\geqslant \widetilde{\angle}_kpqr$ for any $[pq]$ and $[qr]$ in $X$. If
$\angle pqr$ is bad (i.e. $\angle pqr<\widetilde{\angle}_kpqr$), it will be easier to locate an $o\in [qr]\cap  S_p$ via the following lemma
(refer to Subsection 6.1.2 for $S_p$).

\begin{lemma}[{[HSW]}]\label{core1}
Given $[pr_1]$ and $[r_1r_2]$ in $X$, if $\angle pr_1r_2$ is bad, then there is $s\in [r_1r_2]^\circ$ with $|ps|<\max\{|pr_1|,|pr_2|\}$ \footnote{For $k>0$, analogue modifications in footnotes 8 and 9 are required.} such that for any $[ps]$, $\angle psr_1$ or $\angle psr_2$ is bad; in particular, for each $i$, if $|r_is|\leqslant|r_ip|$, then $\angle psr_i$ is bad.
\end{lemma}

To see Lemma 6.12 is easier than to see Lemma 6.9. Note that, for any $[pq],[rs]\subset X$ with $r\in[pq]^\circ$, $\angle prs+\angle qrs=\pi$ (see Subsection 1.3); and for any $[pq],[qr]\subset X$ with $r_i\in[qr]$ and $r_i\to q$ as $i\to\infty$, $|pr_i|\leqslant|pq|-\cos\angle pqr\cdot|qr_i|+ o(|qr_i|)$ (cf. Proposition 6.5).

Furthermore, similar to the proof of Theorem 6.1, starting with the point $o$ we can find $\{o_i\}_{i=1}^\infty$ with $|po_i|$ decreasing as $i\to\infty$.  However, without the completeness, $o_i$ might not have a limit point.
This forces us to explore much deeper analysis on the badness of $\angle pqr$. Indeed, we can do it as follows due to Lemma 3.4; and, by the new analysis, it turns out that we can derive a contradiction
along some minimal geodesic starting from $p$ without involving the points in $S_p$.

\begin{lemma}\label{equiv}
		Let $[pq]\subset X$, and let $\gamma(t)|_{[0,l]}\subset X$ be an arc-length parameterized minimal geodesic with $\gamma(0)=q$.
        If $\angle pq\gamma(l)$ is bad, then there is $[p\gamma(t)]$ with $t\in [0, l)$, a real number $c>0$, and $\tau_i\to 0^+$ as $i\to\infty$  such that
$$|p\gamma(t+\tau_i)|>|p\gamma(t)|-\cos\angle p\gamma(t)\gamma(l) \cdot \tau_i+\frac{1}{2}\ct_k(|p\gamma(t)|)\sin^2\angle p\gamma(t)\gamma(l)\cdot \tau_i^2+c\cdot\tau_i^2.$$
\end{lemma}

\begin{proof} It suffices to show that $\angle pq\gamma(l)\geqslant \widetilde{\angle}_k pq\gamma(l)$ if
$$f_k(|p\gamma(t+\tau)|)\leqslant f_k(|p\gamma(t)|)-\text{sn}_k(|p\gamma(t)|)\cos\angle p\gamma(t)\gamma(l) \cdot \tau+\frac12(1-kf_k(|p\gamma(t)|))\cdot \tau^2+o(\tau^2)\eqno{(6.3)}$$
for all $[p\gamma(t)]$ with $t\in (0,l)$ and sufficiently small $\tau>0$ (cf. (1.2) and (1.3)). As mentioned above, $\angle p\gamma(t)\gamma(0)+\angle p\gamma(t)\gamma(l)=\pi$;
and $|p\gamma(t-\tau)|\leqslant|pq|-\cos\angle p\gamma(t)\gamma(0)\cdot\tau+ o(\tau)$, or equivalently, $f_k(|p\gamma(t-\tau)|)\leqslant f_k(|p\gamma(t)|)-\text{sn}_k(|p\gamma(t)|)\cos\angle p\gamma(t)\gamma(0)\cdot\tau+ o(\tau)$.
These imply that $\liminf\limits_{\tau \to 0^+}\frac{f_k(|p\gamma(t-\tau)|)-f_k(|p\gamma(t)|)}{-\tau}\geqslant-\text{sn}_k(|p\gamma(t)|)\cos\angle p\gamma(t)\gamma(l)$,
which plus (6.3) means that $g_+''(t)+kg(t)\leqslant 1$ in the support sense where $g(t)\triangleq f_k(|p\gamma(t)|)$.

Consider an arc-length parameterized minimal geodesic $\tilde\gamma(t)|_{[0,l]}$ and a point $\tilde p$ in $\mathbb S^2_k$ with
$|\tilde p\tilde\gamma(0)|=|p\gamma(0)|$ and $|\tilde p\tilde\gamma(l)|=|p\gamma(l)|$.
Let $\tilde g(t)\triangleq f_k(|\tilde p\tilde\gamma(t)|)$. Similar to (4.4), $\left(g(t)-\tilde g(t)\right)_+''+k\left(g(t)-\tilde g(t)\right)\leqslant 0$ in the support sense
for all $t\in (0,l)$, and thus by Lemma 3.4 $f_k(|p\gamma(t)|)\geqslant f_k(|\tilde p\tilde\gamma(t)|)$, i.e. $|p\gamma(t)|\geqslant |\tilde p\tilde\gamma(t)|$ (here it needs $l<\frac{\pi}{\sqrt{k}}$ if $k>0$). This implies
$\angle pq\gamma(l)\geqslant \widetilde{\angle}_k pq\gamma(l)$ because $|p\gamma(\tau)|\leqslant|pq|-\cos\angle pq\gamma(l)\cdot\tau+ o(\tau)$ and
$|\tilde p\tilde \gamma(\tau)|=|\tilde p\tilde \gamma(0)|-\cos\widetilde{\angle}_k pq\gamma(l)\cdot\tau+ o(\tau)$ for sufficiently small $\tau>0$.
\end{proof}

Moreover, our proof for Theorem 6.11 needs the following convergence of minimal geodesics which is due to the bounded curvature from below.
Note that by (1.6), for any $x\in X$, it is well defined that
$$R(x)\triangleq\max\{R|\text{ $\angle bac\geqslant\widetilde{\angle}_k bac$ for any $[ab], [ac]\subset B(p,R)$}\}>0,$$
and $R(\cdot)$ is continuous on $X$.

\begin{lemma}\label{conv}
Let $r\in[pq]^\circ\subset X$. If $q_i\to q$ and $r_i\to r$ as $i\to \infty$, then $[r_iq_i]\to [rq]$  for any $[r_iq_i]$.
\end{lemma}

\begin{proof}
It suffices to show that, for any large $N\in \Bbb Z^+$ and $0\leqslant j\leqslant N$,  $r_{ij}\to r_j$ as $i\to\infty$, where $r_{ij}\in [r_iq_i]$ with $|r_ir_{ij}|=\frac{j}{N}|r_iq_i|$
and $r_{j}\in [rq]$ with $|rr_{j}|=\frac{j}{N}|rq|$. Note that by the compactness of $[rq]$, there is $\delta>0$ such that  $R(t)>\delta$ for all $t\in[rq]$. For sufficiently large $N$, we can assume that $\frac{1}{N}|rq|<\delta$,
and there is $r^\prime\in[pr]$ such that $|r^\prime r|=\frac{1}{N}|rq|$. Note that, as $i\to\infty$, $|rr_{i1}|\to |rr_1|=|r'r|$, and $|r^\prime r|+|rr_{i1}|-|r^\prime r_{i1}|\to 0$ (otherwise it will contradict the minimality of $[r'q]$), which means that $\widetilde{\angle}_kr^\prime rr_{i1}\to\pi$.
On the other hand, for large $i$, it is clear that $\{r', r_1, r_{i1}\}\subset B(r, \delta)$, and thus
$$\widetilde{\angle}_kr^\prime rr_{i1}+\widetilde{\angle}_kr_{i1}rr_1\leqslant \angle r^\prime rr_{i1}+\angle r_{i1}rr_1=\pi.\eqno{(6.4)}$$
Plugging `$\widetilde{\angle}_kr^\prime rr_{i1}\to\pi$' into (6.4), we get $\widetilde{\angle}_kr_{i1}rr_1\to 0$ as $i\to\infty$, which implies $r_{i1}\to r_1$.
Then based on this, we can similarly verify that $r_{ij}\to r_j$ for $j\geqslant 2$ as $i\to\infty$ one by one.
\end{proof}

\renewcommand{\proofname}{Proof of Theorem \ref{6.11}}\begin{proof}

As mentioned under Theorem 6.11, one just needs to show $\angle pqr\geqslant \widetilde{\angle}_kpqr$ for any $[pq]$ and $[qr]$ in $X$.
If $\angle pqr<\widetilde{\angle}_kpqr$ (which implies $p\not\in [qr]$), then by Lemma 6.13 there is $[ps]$ with $s\in [qr]\setminus\{r\}$, a real number $c>0$, and $\{r_i\}_{i=1}^\infty\subset[sr]^\circ$ with $r_i\to s$ as $i\to\infty$ such that
$$	|pr_i|>|ps|-\cos\angle psr\cdot|sr_i|+\frac{1}{2}\ct_k(|ps|)\sin^2\angle psr\cdot|sr_i|^2+c\cdot |sr_i|^2.	$$
Then for $p_0\in [ps]^\circ$ sufficiently close to $p$, we have that
$$|p_0r_i|\geqslant|pr_i|-|pp_0|>|p_0s|-\cos\angle p_0sr_i\cdot|sr_i|+\frac{1}{2}\ct_k(|p_0s|)\sin^2\angle p_0sr_i\cdot|sr_i|^2+\frac{c}{2}\cdot |sr_i|^2$$
(note that $\angle p_0sr_i=\angle psr$ and $\ct_k(|p_0s|)>\ct_k(|ps|)$), or equivalently,
$$f_k(|p_0r_i|)>f_k(|p_0s|)-\sn_k(|p_0s|)\cos\angle p_0sr_i\cdot|sr_i|+\frac{1}{2}(1-kf_k(|p_0s|))\cdot|sr_i|^2+\bar c\cdot |sr_i|^2 \eqno{(6.5)}$$
with $\bar c>0$. On the other hand, 		
$$f_k(|p_0r_i|)=f_k(|p_0s|)-\sn_k(|p_0s|)\cos\widetilde\angle_k p_0sr_i\cdot|sr_i|+\frac{1}{2}(1-kf_k(|p_0s|))\cdot|sr_i|^2+o(|sr_i|^2) \eqno{(6.6)}$$
(see (1.3)). From (6.5) and (6.6), it is easy to see that $\liminf\limits_{i\to\infty}\frac{\cos\angle p_0sr_i-\cos \widetilde{\angle}_kp_0sr_i}{|sr_i|}>0$, which implies $\angle p_0sr_i<\widetilde{\angle}_kp_0sr_i$ (i.e. $\angle p_0sr_i$ is bad) for large $i$.
		
By the compactness of $[p_0s]\ (\subset[ps])$, there is $\delta>0$ such that  $R(t)>3\delta$ for all $t\in[p_0s]$. Set a partition of $[sp_0]$: $s=s_0, s_1, \cdots ,s_n=p_0$, where $|s_js_{j+1}|=\delta$ for $0\leqslant j\leqslant n-2$ and $\delta<|s_{n-1}p_0|\leqslant2\delta$.
{\bf Claim}: for each $0\leqslant j\leqslant n-1$, there is $t_j$ and $\bar t_j$ with $|s_jt_j|\ll \delta$ and $|s_j\bar t_j|\ll \delta$ such that $\angle p_0t_j\bar t_j$ is bad.
Note that the badness of $\angle p_0t_{n-1}\bar t_{n-1}$ contradicts `$t_{n-1},\bar t_{n-1}\in B(p_0, R(p_0))$'.

We now just need to verify the claim. When $j=0$, we can let $t_0=s$ and $\bar t_0=r_i$ for any large $i$. We next find out $t_1$ and $\bar t_1$.
Note that any $\triangle t_0s_1\bar t_0$ is a good triangle (i.e. each angle of it is not bad).
Then by Lemma 2.1 (similar to Step II in the proof of Theorem 6.1), $\angle p_0s_1\bar t_0$ is a bad angle; and then by Lemma 6.12 there is $[p_0u_{1i}]$ with $u_{1i}\in [s_1\bar t_0]^\circ$ such that $\angle p_0u_{1i}s_1$ is bad. By Lemma 6.14, $[s_1u_{1i}]\to [s_1u_1]$ with $u_1\in [s_1t_0]$ as $i\to\infty$; and by Lemma 6.14 again $[p_0u_{1i}]\to[p_0u_1]$ \footnote{Note that if $u_1=s$,
it might not be true that $[pu_{1i}]\to [ps]$. This is an essential reason why we need a $p_0$.}.
If $u_1=s_1$, then we can set $t_1=u_{1i}$ for any large $i$ and $\bar t_1=s_1$. If $u_1\neq s_1$, we can let $t_1\in [p_0u_{1i}]$ with $|p_0t_1|=|p_0s_1|$ for any large $i$ and $\bar t_1=s_1$,
because any $\triangle \bar t_1u_{1i}t_1$ ($\subset B(s_1, R(s_1))$) is also a good triangle and so by Lemma 2.1 again $\angle p_0t_1\bar t_1$ is bad.
Then starting from $t_1$ and $\bar t_1$, one by one we can similarly locate the desired $t_j$ and $\bar t_j$ for $2\leqslant j\leqslant n-1$, where
$\bar t_j\in [p_0t_{j-1}]$ with $|p_0\bar t_j|=|p_0s_j|$. (Hint: the requirement of `$\delta<|s_{n-1}p_0|\leqslant2\delta$' is just to guarantee the existence of $u_{(n-1)i}$ with
$\angle p_0u_{(n-1)i}\bar t_{n-1}$ bad by Lemma 6.12).
\end{proof}

%%%%%%%%%%%%%%%%%%%%%%%%%%%%%%%%%%%%%%%%%%%%%%%%%%%%%%%%%%%%%%%

\vskip8mm

\noindent School of Mathematical Sciences (and Lab. Math. Com.
Sys.), Beijing Normal University, Beijing, 100875,
People's Republic of China

\noindent E-mail: 13716456647@163.com; suxiaole$@$bnu.edu.cn; wyusheng$@$bnu.edu.cn


\begin{thebibliography}{MMM}

\bibitem[AKP]{AKP} S. Alexander, V. Kapovitch and A. Petrunin, {\it Alexandrov Geometry}, 2010.

\bibitem[BGP]{BGP}
Yu. Burago, M. Gromov and G. Perel$'$man, {\it A.D. Aleksandrov spaces
with curvature bounded below}, Uspeckhi Mat. Nank, 47(2): 3-51, 1992.

\bibitem[HSW]{HSW} S. Hu, X. Su and Y. Wang, {\it A proof of Toponogov's theorem in Alexandrov geometry}, Proc. of AMS, 151(4): 1743-1748, 2023.

\bibitem[MV]{MV}
A.D. Miller and R. Vyborny, {\it Some remarks on functions with one-sided derivatives}, The American Mathematical Monthly, 93(6): 471-475, 1986.

\bibitem[Na]{Na} I.P. Natanson, {\it Theory of funcions of a real variable}, Teoria functsiy veshchestvennoy peremennoy. Frederick Ungar Pub. Co, 1960.

\bibitem[Pe]{Pe}
G. Perel$'$man, {\it A.D. Alexandrov spaces with curvature bounded below II}, 1992, preprint.

\bibitem[PP]{PP}
G. Perel$'$man and A. Petrunin, {\it Quasigeodesics and Gradient Curves in Alexandrov spaces}, 1994,
www.math.psu.edu/Petrunin/papers/.

\bibitem[Pet]{Pet} P. Petersen, {\it Riemannian Geometry} (2nd Edition), Springer, United States, 2006.

\bibitem[Petr1]{Petr1} A. Petrunin, {\it Quasigeodesics in multidimensional Alexandrov spaces}, Thesis (Ph.D.)-University of Illinois at Urbana-Champaign, 1995.

\bibitem[Petr2]{Petr2} A. Petrunin, {\it A globalization for non-complete but geodesic spaces}, Math. Ann. 366: 387-393, 2016.

\bibitem[Pl]{Pl} C. Plaut, {\it Spaces of Wald-Berestovskii curvature bounded below}, J. Geom. Anal. 6: 113-134, 1996.

\bibitem[Sh]{Sh} K. Shiohama, {\it An introduction to the geometry of Alexandrov spaces}, 1992.


\end{thebibliography}
\end{document}